\documentclass{article}[11pt]

\usepackage{amsmath} 
\usepackage{amssymb}
\usepackage{amsthm}
\usepackage{accents}
\usepackage{xcolor}
\usepackage{tikz-cd}

\newtheorem{Th}{Theorem}[section] 
\newtheorem{Prop}{Proposition}[section]   
\newtheorem{Lem}{Lemma}[section]   
\newtheorem{Coro}{Corollary}[section]   
\newtheorem{Def}{Definition}  
\newtheorem{Rem}{Remark}[section]

\newcommand{\R}{\mathbb{R}}
\newcommand{\Z}{\mathbb{Z}}
\newcommand{\C}{\mathbb{C}}

\newcommand{\N}{\mathbb{N}}
\newcommand{\Sz}{\mathcal{S}}
\newcommand{\LL}{\mathcal{L}}
\newcommand{\D}{\mathcal{D}}
\newcommand{\E}{\mathcal{E}}

\newcommand{\U}{\mathcal{U}}

\newcommand{\aW}{{\accentset{\,\circ}W}}

\newcommand{\x}{\langle x\rangle}

\newcommand{\p}{{\rm p}}

\newcommand{\id}{\text{\rm id}}

\newcommand{\Div }{\mathop{\rm div}\nolimits}

\newcommand{\dt}[1]{\accentset{\mbox{\bfseries .}}{#1}}
\newcommand{\bfdot}{{\mbox{\bfseries .}}}
\newcommand{\dd}{{\rm d}}

\newcommand{\3}{\vert}

\begin{document}

\title{A holomorphic Lie-Trotter product formula for nonlinear local semigroups}   
 
\author{P. Topalov} 

\maketitle

\begin{abstract}  
We prove a new, holomorphic version of the Lie-Trotter product formula for nonlinear local semigroups with a parameter.
We show that the product converges to the solutions of a nonlinear evolutions equation with a parameter. 
The solution is holomorphic in time on a centered at zero conic set in the complex plane and depends continuously 
on the parameter, and analytically on the initial data.
\end{abstract}   


\section{Introduction}\label{sec:introduction}
Let $U$ be a non-empty open set in a complex Banach space $(X,\|\cdot\|)$, $0<\vartheta_0<\pi/2$, and $0<T\le\infty$. 
Consider the conic set in $\C$,
\begin{equation}\label{eq:the_cone}
\mathbb{S}_{\vartheta_0}:=\big\{z\in\C\setminus\{0\}\,\big|\,|\arg z|<\vartheta_0\big\},
\end{equation}
as well as the sets
\begin{equation}\label{eq:other_cones}
\mathbb{S}_{\vartheta_0,T}:=\big\{z\in\mathbb{S}_{\vartheta_0}\,\big|\,|z|<T\big\},\quad
\mathbb{S}_{\vartheta_0,0}:=\mathbb{S}_{\vartheta_0}\sqcup\{0\},\quad\text{\rm and}\quad
\mathbb{S}_{\vartheta_0,0,T}:=\mathbb{S}_{\vartheta_0,T}\sqcup\{0\} .
\end{equation}
We say that a continuous map
\[
S : \mathbb{S}_{\vartheta_0,0,T}\times U\to X
\] 
is a {\em local semigroup} if it satisfies the following {\em local semigroup property in $U$}: 
for any $x\in U$ we have that $S_t(x)|_{t=0}=x$ and $S_t\circ S_s(x)=S_{t+s}(x)$ provided that 
$S_s(x)\in U$ and $s,t,s+t\in\mathbb{S}_{\vartheta_0,0,T}$. 
Here, and in what follows, we set $S_t(x):=S(t,x)$.
Let $\widetilde{X}$ be a complex Banach space such that $X\subseteq\widetilde{X}$, and the inclusion is bounded.
We are interested in the equation
\begin{equation}\label{eq:pde_introduction}
\frac{\dd u}{\dd t}=\E(u)+\nu\,\Xi(u),\quad u|_{t=0}=x\in U,
\end{equation}
where $0\le\nu\le 1$ is a parameter and the maps
\[
\E : U\to X\quad\text{\rm and}\quad\Xi : U\to\widetilde{X}
\]
are analytic. We will assume that $\Xi$ generates a local semigroup 
\begin{equation}\label{eq:S_introduction}
S : \mathbb{S}_{\vartheta_0,0,T}\times U\to X
\end{equation}
such that the map $S : \mathbb{S}_{\vartheta_0,0,T}\times U_Z\to\widetilde{X}$ is differentiable and
\begin{equation*}
\frac{\dd}{\dd t}\Big|_{t=0} S_t(z)=\Xi(z)
\end{equation*}
for any $z\in U\cap Z$, where $Z\subseteq X$ is a densely and boundedly embedded in $X$ Banach space, and
$U_Z=U\cap Z$ is an open set in $Z$.
The derivative at $t=0$ is defined as the limit of $(S_t(z)-z)/t$ in $\widetilde{X}$ as $t\to 0$ 
with $t\in\mathbb{S}_{\vartheta_0,T}$.
(The properties of such maps, called {\em maps on conic sets}, are studied in Appendix \ref{sec:C^k_along_curves}.)
We prove in Proposition \ref{prop:product_formula} and Theorem \ref{th:product_formula} that,
under certain conditions on the local semigroup \eqref{eq:S_introduction} and the local semigroup
\begin{equation}\label{eq:G_introduction}
G : \mathbb{S}_{\vartheta_0,0,T}\times U\to X
\end{equation}
generated by the analytic vector field $\E : U\to X$, 
we have that for any value of the parameter $0\le\nu\le 1$
the equation \eqref{eq:pde_introduction} has a unique solution
$u=u^{(\nu)}\in C\big([0,\tau),X\big)\cap C^1\big([0,\tau),\widetilde{X}\big)$
where the time of existence $0<\tau\le T$ can be chosen locally uniform on the choice
of the initial data  $x\in U$. The solution has strong regularity properties: It extends to
a holomorphic map $u^{(\nu)} : \mathbb{S}_{\vartheta_0,\tau}\to X$, 
$u^{(\nu)}\in C\big(\mathbb{S}_{\vartheta_0,0,\tau},X\big)$, that satisfies 
\eqref{eq:pde_introduction} with one-sided derivative at $t=0$ (cf. Appendix \ref{sec:C^k_along_curves}).
Moreover, the corresponding (locally defined) flow map
\begin{equation}\label{eq:the_flow_map}
\mathbb{S}_{\vartheta_0,0,\tau}\times O\times[0,1]\to X,\quad
(t,x;\nu)\mapsto u^{(\nu)}(t,x),
\end{equation}
where $O\subseteq U$ is an open subset in $U$, is continuous, 
and for any given value of the parameter $0\le\nu\le 1$, the restricted to 
$\mathbb{S}_{\vartheta_0,\tau}\times O$ map
\begin{equation}\label{eq:u^nu}
u^{(\nu)} : \mathbb{S}_{\vartheta_0,\tau}\times O\to X
\end{equation}
is analytic. In particular, for any given $x\in O$ we have that
\begin{equation}\label{eq:viscosity_limit}
u^{(\nu)}\to u^{(0)}\quad\text{\rm as}\quad\nu\to 0+
\end{equation}
where $u^{(0)}\in C\big([0,\tau),X\big)\cap C^1\big([0,\tau),\widetilde{X}\big)$ is the solution
of \eqref{eq:pde_introduction} with $\nu=0$, and the limit is taken in 
$C\big([0,\tau),X\big)\cap C^1\big([0,\tau),\widetilde{X}\big)$.
The convergence in \eqref{eq:viscosity_limit} is, in fact, stronger 
(see Proposition \ref{prop:product_formula} and Corollary \ref{coro:nu}). 
The limit \eqref{eq:viscosity_limit} is an analog of the ``zero viscosity'' limit in fluid dynamics. 
Another important ``computational'' property of the solution \eqref{eq:u^nu} is that 
for any $t\in\mathbb{S}_{\vartheta_0,0,\tau}$, $x\in O$, and $0\le\nu\le 1$, we have that 
\begin{equation}\label{eq:product_formula_limit}
u^{(\nu)}(t,x)=\lim_{n\to\infty}\big(S^{\nu}_{t/n}\circ G_{t/n}\big)^n(x),\quad 
S^{(\nu)}_t(x)=S_{\nu t}(x),
\end{equation}
where the sequence converges in $X$, locally uniformly in $x$ in $O$ and
$t\in\mathbb{S}_{\vartheta_0,0,\tau}$, and uniformly in $0\le\nu\le 1$. 
We refer to Proposition \ref{prop:product_formula} for the details.

\medskip

The paper is organized as follows. In Section \ref{sec:Lie-Trotter} we formulate
our main results -- Proposition \ref{prop:product_formula} and Theorem \ref{th:product_formula}.
Let us explain the relation between these two statements.
In Proposition \ref{prop:product_formula} we prove that, under certain conditions on the local semigroups
\eqref{eq:S_introduction} and \eqref{eq:G_introduction}, the equation \eqref{eq:pde_introduction} has
a unique solution that satisfies the properties discussed above. However, one of the conditions
(item (b)) in Proposition \ref{prop:product_formula} requires that each of the generated local semigroups is 
locally equivalent to a {\em quasi-contracion} (Definition \ref{def:quasi-contraction}). 
The quasi-contracion condition is not generic and could be difficult to verify. 
Theorem \ref{th:product_formula} reformulates condition (b) in terms of global in time estimates
on the local semigroup \eqref{eq:S_introduction} (cf. estimates \eqref{eq:dS_t-exponential_growth1} and 
\eqref{eq:dS_t-exponential_growth2}). Such global estimates can be verified in a number of situations, 
where Theorem \ref{th:product_formula} can be applied. 
Proposition \ref{prop:product_formula} and Theorem \ref{th:product_formula} are proved in
Section \ref{sec:Lie-Trotter_proofs}. We start the Section with the proof of two technical Lemmas --
Lemma \ref{lem:equivalent_distance} and Lemma \ref{lem:stability}. 
Lemma \ref{lem:equivalent_distance} shows that the global in time estimates 
\eqref{eq:dS_t-exponential_growth1} and \eqref{eq:dS_t-exponential_growth2}
together with the properties of the generator $\E$ imply that the quasi-contraction condition (b) of 
Proposition \ref{prop:product_formula} holds. In Lemma \ref{lem:stability} we prove that
the successive iterations of the local semigroups \eqref{eq:S_introduction} and \eqref{eq:G_introduction}
is locally bounded, uniformly on the choice of the initial point $x\in U$ and $t\in\mathbb{S}_{\vartheta_0,0,\tau}$
for some $0<\tau\le T$ (cf. the product formula \eqref{eq:stability}). In Lemma \ref{lem:fundamental} we 
prove that if the initial point $x\in U\cap Z$ of the iteration belongs to a dense subset $U\cap Z$ in $U$, 
then the successive iterations converge in $X$, locally uniformly on the choice of the initial point $x\in U\cap Z$ and 
$t\in\mathbb{S}_{\vartheta_0,0,\tau}$ for some $0<\tau'\le\tau$. This can be done uniformly
in the parameter $0\le\nu\le 1$. A crucial observation in the proof of Proposition \ref{prop:product_formula} is 
that Lemma \ref{lem:fundamental} and the fact that the products \eqref{eq:stability_Z} are holomorphic 
functions of $t$ on the conic set $\mathbb{S}_{\vartheta_0,\tau}$ allow us to apply Vitali's theorem 
(with parameter) to conclude that the limit \eqref{eq:product_formula_limit} exists for $x\in U\cap Z$ and
$t$ in the larger cone $\mathbb{S}_{\vartheta_0,\tau}$, uniformly in the parameter $0\le\nu\le 1$. 
We then use the quasi-contraction property and approximation arguments to show that 
the limit also exists for $x\in U$, and that the limit function $u^{(\nu)}$ is a solution of  
\eqref{eq:pde_introduction}. The analytic dependence on the initial data and the remaining 
regularity properties of the constructed solution then follow relatively easily by construction.
In Appendix \ref{sec:C^k_along_curves} we define the notion of $C^k$ maps on conic sets and prove
their main properties. The framework is needed for the systematic treatment of the composition and 
the regularity at $t=0$ of maps defined on conic sets.

\medskip

Different versions of the Lie-Trotter product formula play an important role in mathematics and its applications.
Here we only mention the seminal work of Ebin and Marsden \cite{EM}, where a product formula argument was 
used for proving the local well-posedness of the Navier-Stokes equation on closed manifolds. 
In mathematical physics, the product formula is used for the justification of the Feynman path integral 
(cf., e.g., \cite{Nelson,Takhtajan}).
Proposition \ref{prop:product_formula} and Theorem \ref{th:product_formula} can be considered
as analogs of \cite[Theorem 13.1]{EM} and \cite[Theorem 2.1 and Theorem 5.1]{Marsden}, extended
to the case of holomorphic nonlinear local semigroups. In addition to the extention of these results to
holomorphic nonlinear local semigroups, our main motivation is to obtain a version of 
\cite[Theorem 13.1]{EM} that does not require that the nonlinear local semigroup generated 
by $\Xi$ is a contraction. The reason is that in a number of applications of these results to fluid equations, 
the generator $\Xi$ is the Laplace operator conjugated by a diffeomorphism of $\R^d$, $d\ge 2$, 
which has a specified growth/decay condition at infinity \cite{McOwenTopalov3,McOwenTopalov5} 
(see also \cite{McOwenTopalov2,McOwenTopalov4,Top1}, as well as 
the classical works of Arnold \cite{Arnold}, G\"unther \cite{Gunther}, and Wolibner \cite{Wol}).
Such diffeomorphisms are modeled on Sobolev spaces with weights, and the generated 
by the conjugated Laplace operator heat flow on such spaces is not a contraction,
but a quasi-contraction of the form $\|S_t(x)\|\le M e^{\beta|t|}$, $x\in X$, with a constant $M$
not necessarily equal to one (see \cite[\S 2, \S 3]{McOwenTopalov5}).
The assumptions in \cite[Theorem 2.1 and Theorem 5.1]{Marsden} also do not hold in this case.
(This concerns condition (iii) in \cite[Theorem 2.1]{Marsden}, which holds on unbounded sets,
and the requirement in \cite[Theorem 5.1]{Marsden} that the second derivatives in 
the spatial direction of the nonlinear local semigroup $S_t$ are uniformly 
bounded on bounded sets.) In addition, in the majority of the cases, the existence of 
the limit \eqref{eq:viscosity_limit} is proved under the assumption that the right side of 
\eqref{eq:pde_introduction} generates a local semigroup (cf., e.g., \cite[Theorem  9]{Nelson}, 
\cite[\S 3.4.2]{MB}, \cite{Trotter}). 
Finally, note that there is a well-developed theory of nonlinear semigroups 
(see, e.g., \cite{Komura,Kato0,BrezisPazy,LL,Oharu1,HaydenMassey} and the references therein). 
However, these fundamental results do not consider dependence on parameters and usually
assume that the generators are monotone nonlinear operators, or other analogous conditions, 
which in our situation are not satisfied.
For the zero viscosity limit for fluid equations, we refer 
to \cite{Kato1,Kato2,Swann,McGrath,MB,Pazy}.

\medskip

In the case when $X$ and $Y$ are Banach spaces over $\C$, we denote by $\LL(X,Y)$ the space
of complex linear maps equipped with the uniform norm. For a given metric space $\Omega$ we
denote by $C_b\big(\Omega,X\big)$ the Banach space of bounded continuous maps $\Omega\to X$ 
equipped with the supremum norm.

\section{The holomorphic Lie-Trotter product formula}\label{sec:Lie-Trotter}
In this Section we state a holomorphic version of the Lie-Trotter product formula for nonlinear local semigroups.
Let $U$ be a non-empty open set in a complex Banach space $(X,\|\cdot\|)$, $0<\vartheta_0<\pi/2$, and $0<T\le\infty$,
and $N$ be an non-empty open set in $X$. A distance function $d : N\times N\to\R_{\ge 0}$ on $N$ is 
{\em equivalent to the standard distance} on $N$ if there exist constants $0<L_1<L_2<\infty$ such that
\begin{equation}\label{eq:equivalent_distances}
L_1\,d(x,y)\le\|x-y\|\le L_2\,d(x,y)
\end{equation}
for any $x,y\in N$.

\begin{Def}\label{def:quasi-contraction}
Assume that $S : \mathbb{S}_{\vartheta_0,0,T}\times U\to X$ is a local semigroup.
Then, we say that $S$ is {\em locally equivalent to a quasi-contraction in $U$} if 
for any $x_0\in U$ there exist an open neighborhood $N(x_0)\subseteq U$ of $x_0$ in $X$,
a distance function $d : N(x_0)\times N(x_0)\to\R_{\ge 0}$ on $N(x_0)$ equivalent to the standard distance on $N(x_0)$, 
and constants $\beta>0$ and $0<T_0\le T$, such that if the open set $V\subseteq N(x_0)$ in $X$ and $0<\tau\le T_0$ are 
such that $S_t(V)\subseteq N(x_0)$ for any $t\in\mathbb{S}_{\vartheta_0,0,\tau}$ then 
\begin{equation*}
d\big(S_t(x),S_t(y)\big)\le e^{\beta|t|}d(x,y)
\end{equation*}
holds for any $x,y\in V$ and $t\in\mathbb{S}_{\vartheta_0,0,\tau}$. The map $S$ is {
\em globally equivalent to a quasi-contraction in $U$} if the above holds for $N(x_0)$ replaced by $U$.
\end{Def}

\noindent  The maps $S : \mathbb{S}_{\vartheta_0,0,T}\times U\to X$ and 
$G : \mathbb{S}_{\vartheta_0,0,T}\times U\to X$ are called 
{\em locally equivalent to a quasi-contraction in $U$ with the same distance function} if 
for any $x_0\in U$ there exist an open neighborhood $N(x_0)$, a distance function 
$d : N(x_0)\times N(x_0)\to\R_{\ge 0}$ as in Definition \ref{def:quasi-contraction}, and 
constants $\beta>0$ and $0<T_0\le T$, so that if the open set $V\subseteq N(x_0)$ in $X$ and $0<\tau\le T_0$ are
such that $G_t(V),S_t(V)\subseteq N(x_0)$ for any $t\in\mathbb{S}_{\vartheta_0,0,\tau}$ then 
\begin{equation}\label{eq:quasi-contraction}
d\big(G_t(x),G_t(y)\big)\le e^{\beta |t|}d(x,y)\quad\text{\rm and}\quad
d\big(S_t(x),S_t(y)\big)\le e^{\beta |t|}d(x,y)
\end{equation}
for any $x,y\in V$ and $t\in\mathbb{S}_{\vartheta_0,0,\tau}$.
If \eqref{eq:quasi-contraction} holds for $N(x_0)$ replaced by $U$ then we say that $G$ and $S$ are
{\em globally equivalent to a quiasi-contraction in $U$ with the same distance function.} 

Let us formulate our first result.
For the convenience of the reader, the various maps appearing in  Proposition \ref{prop:product_formula} below 
are summarized in the commutative diagram \eqref{eq:diagram}.
In the diagram, $t\in\mathbb{S}_{\vartheta_0,0,T}$ is fixed, the (hooked) vertical arrows 
denote embeddings of subsets of Banach spaces, and the various restrictions of a given map are 
denoted by the same letter.
For an open set $U\subseteq X$ and a boundedly embedded Banach space $Z\subseteq X$ 
we denote by $U_Z\subseteq X$ the open set $U\cap Z$ in $Z$. (Note that $U_Z$ need not be bounded in $Z$
when $U$ is bounded in $X$.)

\begin{Prop}\label{prop:product_formula}
Let $(Z,\|\cdot\|_Z)$ and $(X,\|\cdot\|)$ be complex Banach spaces such that the inclusion $Z\subseteq X$
is bounded and $Z$ is dense in $X$. 
Assume that there exist a non-empty open set $U$ in $X$ and a constant $T>0$ such that the maps
\begin{equation}\label{eq:G}
G : \mathbb{S}_{\vartheta_0,0,T}\times U\to X\quad\text{\em and}
\quad G : \mathbb{S}_{\vartheta_0,0,T}\times U_Z\to Z,
\end{equation}
and
\begin{equation}\label{eq:S}
\quad S : \mathbb{S}_{\vartheta_0,0,T}\times U\to X\quad\text{\em and}
\quad S : \mathbb{S}_{\vartheta_0,0,T}\times U_Z\to Z
\end{equation}
are local semigroups, and the restricted maps
\begin{equation}\label{eq:S,G-restricted}
G : \mathbb{S}_{\vartheta_0,T}\times U\to X\quad\text{\em and}\quad
S : \mathbb{S}_{\vartheta_0,T}\times U\to X
\end{equation}
are analytic.
In addition, we assume:
\begin{itemize}
\item[$(a)$] {\sc Regularity at $t=0$:} 
$1)$ The maps in \eqref{eq:G} are $C^2$ on conic sets (see Definition \ref{def:C^k_along_curves_complex} in 
Appendix \ref{sec:C^k_along_curves}); 
$2)$ The map $S : \mathbb{S}_{\vartheta_0,0,T}\times U_Z\to X$ is $C^2$ on conic sets and the map
$S : \mathbb{S}_{\vartheta_0,0,T}\times U_D\to Z$, 
where the Banach space $D\subseteq Z$ is densely and boundedly embedded in $Z$,
is $C^1$ on conic sets.

\begin{equation}\label{eq:diagram}
\begin{tikzcd}
\widetilde{X}&U\arrow[r, "S_t{,}\,G_t"]\arrow[l, swap, "\Xi{,}\,\,\E"]&X\\
&U_Z\arrow[u, hook]\arrow[r, "S_t{,}\,G_t"]\arrow[ur,"S_t"]&Z\arrow[u, hook]\\
&U_D\arrow[u, hook]\arrow[ur, swap, "S_t"]                 
\end{tikzcd}
\end{equation}

\item[$(b)$] {\sc Quasi-contraction:} 
$1)$ The maps $G : \mathbb{S}_{\vartheta_0,0,T}\times U\to X$ and 
$S : \mathbb{S}_{\vartheta_0,0,T}\times U\to X$ are locally equivalent to a quasi-contraction in $U$ with 
the same distance function; $2)$ The maps $G : \mathbb{S}_{\vartheta_0,0,T}\times U_Z\to Z$ and 
$S : \mathbb{S}_{\vartheta_0,0,T}\times U_Z\to Z$ are locally equivalent to a quasi-contraction in $U_Z$ with 
the same distance function.

\item[$(c)$] {\sc Generators:} 
For a given $x,z\in U_Z$ consider the generators\footnote{See \eqref{eq:derivative_at_t=0} in 
Appendix \ref{sec:C^k_along_curves} for the definition of the derivative at $t=0$.}
\begin{equation*}
\E(x):=\frac{\dd}{\dd t}\Big|_{t=0}G_t(x)\quad\text{\rm and}\quad\Xi(z):=\frac{\dd}{\dd t}\Big|_{t=0} S_t(z)
\end{equation*}
of the local semigroups \eqref{eq:G}, \eqref{eq:S}, and assume that  $\Xi : U_Z\to X$ extends to 
an analytic map
\begin{equation}\label{eq:Xi-regularity}
\Xi : U\to\widetilde{X}
\end{equation}
where $\widetilde{X}$ are Banach spaces such that the embedding $X\subseteq\widetilde{X}$ is bounded.
(The analyticity of $\E : U\to X$ follows since the first map in \eqref{eq:G} is $C^2$ on conic sets 
(cf. Proposition \ref{prop:analytic_derivatives} in Appendix \ref{sec:C^k_along_curves}).)
\end{itemize}
Then, for any given $x_0\in U$ there exist $0<\tau\le T$ and an open neighborhood $O(x_0)\subseteq U$ of $x_0$ in $X$
such that for any given $x\in O(x_0)$ and $0\le\nu\le 1$ the equation
\begin{equation}\label{eq:pde}
\frac{\dd u}{\dd t}=\E(u)+\nu\,\Xi(u),\quad u|_{t=0}=x,
\end{equation}
has a unique solution
\begin{equation}\label{eq:solution_real}
u=u^{(\nu)}\in C_b\big([0,\tau),X\big)\cap C^1_b\big([0,\tau),\widetilde{X}\big)
\end{equation}
that depends Lipschitz continuously on the initial data $x\in O(z_0)$.
The solution \eqref{eq:solution_real} extends to a holomorphic map
$u^{(\nu)} : \mathbb{S}_{\vartheta_0,\tau}\to X$, $u^{(\nu)}\in C_b\big(\mathbb{S}_{\vartheta_0,0,\tau},X\big)$,
that satisfies equation \eqref{eq:pde} with one-sided derivative of $u^{(\nu)}$ at $t=0$ 
(cf. \eqref{eq:derivative_at_t=0} in Appendix \ref{sec:C^k_along_curves}). The flow map
\begin{equation*}
H^{(\nu)} : \mathbb{S}_{\vartheta_0,0,\tau}\times O(x_0)\to X,\quad(t,x)\mapsto u^{(\nu)}(t,x),
\end{equation*} 
is continuous and analytic when restricted to $\mathbb{S}_{\vartheta_0,\tau}\times O(x_0)$.
Moreover, $H^{(\nu)}$ satisfies the local semigroup property in $O(x_0)$, the map
\begin{equation}\label{eq:H-continuous}
H : \mathbb{S}_{\vartheta_0,0,\tau}\times O(x_0)\times[0,1]\to X,\quad(t,x;\nu)\mapsto u^{(\nu)}(t,x),
\end{equation}
is continuous, and there exists a constant $\beta>0$ such that
\begin{equation}\label{eq:H-Lipschitz*}
d\big(H_t^{(\nu)}(x_1),H_t^{(\nu)}(x_2)\big)\le e^{2\beta|t|}d(x_1,x_2)
\end{equation}
for any $x_1,x_2\in O(x_0)$ and $t\in\mathbb{S}_{\vartheta_0,0,\tau}$, where $d : O(x_0)\times O(x_0)\to\R_{\ge 0}$ is 
a distance function equivalent to the induced by $X$ distance on $O(x_0)$.
For any $t\in\mathbb{S}_{\vartheta_0,0,\tau}$ and $x\in O(x_0)$ 
we have that $(S_{t/n}\circ G_{t/n}\big)^n(x)\in U$ for any $n\ge 1$ and
\begin{equation}\label{eq:Lie-Trotter}
H_t^{(\nu)}(x)=\lim_{n\to\infty}\big(S_{t/n}^{(\nu)}\circ G_{t/n}\big)^n(x),\quad S_t^{(\nu)}(x):=S_{\nu t}(x),
\end{equation}
where the limit is taken in $X$ and the convergence is uniform in $x$ on compact subsets of $O(x_0)$,
locally uniform in $t\in\mathbb{S}_{\vartheta_0,0,\tau}$, and uniform in $1\le\nu\le 1$.
\end{Prop}

The Corollary below follows easily from the continuity of \eqref{eq:H-continuous} and Cauchy's estimate.
Recall that, by definition, $f_n\to f$ as $n\to\infty$ in $C^k\big((0,\tau),\widetilde{X}\big)$ if 
$f_n\to f$ as $n\to\infty$ in $C^k\big(J,\widetilde{X}\big)$ for any compact interval $J\subseteq(0,\tau)$.

\begin{Coro}\label{coro:nu}
Let $u^{(\nu)}\in C_b\big([0,\tau),X\big)\cap C^1_b\big([0,\tau),\widetilde{X}\big)$ be the unique solution of  \eqref{eq:pde}
constructed in Proposition \ref{prop:product_formula}.
Then, we have that $u^{(\nu)}\to u^{(0)}$ as $\nu\to 0+$
in $C_b\big( [0,\tau),X\big)\cap C^1_b\big( [0,\tau),\widetilde{X}\big)$.
Moreover,
$u^{(\nu)}\in C_b\big([0,\tau),X\big)\cap C^1_b\big([0,\tau),\widetilde{X}\big)\cap C^\infty\big((0,\tau),\widetilde{X}\big)$
and for any $k\ge 1$ we have that  $u^{(\nu)}\to u^{(0)}$ as $\nu\to 0+$ in
$C_b\big([0,\tau),X\big)\cap C^1_b\big([0,\tau),\widetilde{X}\big)\cap C^k\big((0,\tau),\widetilde{X}\big)$.
\end{Coro}

The quasi-contraction condition $(b)$ in Proposition \ref{prop:product_formula} is very often difficult to verify. 
The following Theorem gives a useful sufficient condition that involves global time estimates on the first and 
the second derivatives of $S_t : U\to X$. 

\begin{Th}\label{th:product_formula}
The statement of Proposition \ref{prop:product_formula} holds with the condition $(b)$ in Proposition \ref{prop:product_formula}
replaced by the following two conditions: 
$1')$ For any given $x_0\in U$ there exists an open neighborhood $N_X(x_0)$ of $x_0$ in $X$ such that
the map $S : \mathbb{S}_{\vartheta_0}\times N_X(x_0)\to X$ is a local semigroup,
for any given $t\in\mathbb{S}_{\vartheta_0}$ the map $S_t : N_X(x_0)\to X$ is analytic, 
$S_t\big(N_X(x_0)\big)\subseteq N_X(x_0)$, and there exist constants $M_X,\beta_X>0$ such that
\begin{equation}\label{eq:dS_t-exponential_growth1}
\|\dd_xS_t\|_{\LL(X,X)}\le M_X e^{\beta_X|t|}\quad\text{\rm and}\quad
\|\dd^2_xS_t\|_{\LL(X\times X,X)}\le M_X e^{\beta_X|t|}
\end{equation}
hold for any $x\in N_X(x_0)$ and for $t\in\mathbb{S}_{\vartheta_0}$;
$2')$ For any given $z_0\in U_Z$ there exists an open 
neighborhood $N_Z(z_0)$ of $z_0$ in $Z$ such that the map
$S : \mathbb{S}_{\vartheta_0}\times N_Z(z_0)\to Z$ is a local semigroup,
for any given $t\in\mathbb{S}_{\vartheta_0}$ the map $S_t : N_Z(z_0)\to Z$ is analytic,
$S_t\big(N_Z(z_0)\big)\subseteq N_Z(z_0)$, and  there exist constants $M_Z,\beta_Z>0$ such that
\begin{equation}\label{eq:dS_t-exponential_growth2}
\|\dd_zS_t\|_{\LL(Z,Z)}\le M_Z e^{\beta_Z|t|}\quad\text{\rm and}\quad
\|\dd^2_zS_t\|_{\LL(Z\times Z,Z)}\le M_Z e^{\beta_Z|t|}
\end{equation}
hold for any $z\in N_Z(z_0)$ and for $t\in\mathbb{S}_{\vartheta_0}$.
\end{Th}

\begin{Rem}
Note that in Theorem \ref{th:product_formula} we do {\em not} assume  that $N_X(x_0)\subseteq U$ or 
that $N_Z(x_0)\subseteq U_Z$. In the applications, the neighborhood $N_X(x_0)$ (resp. $N_Z(z_0)$) is larger 
than $U$ (resp. $U_Z$). In this way, in Theorem \ref{th:product_formula} we assume that the map
\eqref{eq:S} extends to a larger domain $\mathbb{S}_{\vartheta_0}\times N_X(x_0)$ 
(resp. $\mathbb{S}_{\vartheta_0}\times N_Z(z_0)$).
\end{Rem}

\section{Proof of Proposition \ref{prop:product_formula} and 
Theorem \ref{th:product_formula}}\label{sec:Lie-Trotter_proofs} 
For the proof of Theorem \ref{th:product_formula} we need the following Lemma.

\begin{Lem}\label{lem:equivalent_distance}
Let $N_1\subseteq N\subseteq X$ be open sets in a complex Banach space $(X,\|\cdot\|)$. 
Assume that for some $T>0$ and $0<\vartheta_0<\pi/2$ the maps $G : \mathbb{S}_{\vartheta_0,0,T}\times N_1\to X$ and 
$S : \mathbb{S}_{\vartheta_0,0}\times N\to X$ are continuous and satisfy the local semigroup property in $N_1$ and
$N$, respectively. Assume in addition that:
\begin{itemize}
\item[$(i)$] For any given $t\in\mathbb{S}_{\vartheta_0}$ the map $S_t : N\to X$ is analytic,
$S_t(N)\subseteq N$, and there exist constants $M,\beta>0$ such that
\begin{equation}\label{eq:dS_t-exponential_growth}
\|\dd_xS_t\|_{\LL(X,X)}\le M e^{\beta |t|}\quad\text{\rm and}\quad
\|\dd^2_xS_t\|_{\LL(X\times X,X)}\le M e^{\beta |t|}
\end{equation}
for any $x\in N$ and $t\in\mathbb{S}_{\vartheta_0}$.
\item[$(ii)$] The map $G : \mathbb{S}_{\vartheta_0,0,T}\times N_1\to X$ is $C^2$ 
on conic sets\footnote{This implies that $\E : N_1\to X$ is analytic.}
and there exists a constant $L>0$ such that
\begin{equation}\label{eq:E-estimates}
\|\E(x)\|\le L\quad\text{\rm and}\quad\|\dd_x\E\|_{\LL(X,X)}\le L
\end{equation}
for any $x\in N_1$.
\end{itemize}
Then, $S$ and $G$ are locally equivalent to a quasi-contraction in $N_1$ with the same distance function. 
More specifically, for any $x_0\in N_1$ there exist an open neighborhood $U(x_0)\subseteq N_1$ of $x_0$ in $X$,
a constant $0<T_0\le T$, and a distance function $d : U(x_0)\times U(x_0)\to\R_{\ge 0}$ equivalent to 
the standard distance on $U(x_0)$ such that if the open set $V\subseteq U(x_0)$ in $X$ and $0<\tau\le T_0$ are
chosen so that $S_t(V), G_t(V)\subseteq U(x_0)$ for any $t\in\mathbb{S}_{\vartheta_0,0,\tau}$ then 
\begin{equation}\label{eq:S,G-quasi-contraction}
d\big(S_t(x),S_t(y)\big)\le e^{\beta |t|}d(x,y)\quad\text{\rm and}\quad
d\big(G_t(x),G_t(y)\big)\le e^{2ML |t|}d(x,y)
\end{equation}
for any $x,y\in V$ and $t\in\mathbb{S}_{\vartheta_0,0,\tau}$.
\end{Lem}

\proof[Proof of Lemma \ref{lem:equivalent_distance}]
For any given $x\in N$ and $\xi\in X$ define the norm\footnote{See \cite[Theorem 4.1]{Marsden} for the case of $t$ real.} on the tangent space $T_xN=X$  
\begin{equation}\label{eq:Finsler_norm}
|\xi|_x:=\sup\limits_{t\in\mathbb{S}_{\vartheta_0,0}}\big\|e^{-\beta|t|}(\dd_x S_t)(\xi)\big\|\,.
\end{equation}
Without loss of generality, we assume that $M\ge 1$.
It then follows from \eqref{eq:Finsler_norm} and the first estimate in \eqref{eq:dS_t-exponential_growth} that 
for any $x\in N$ and $\xi\in T_xN$ we have that
\begin{equation}\label{eq:Finsler-Euclidean_norms}
\|\xi\|\le|\xi|_x\le  M\|\xi\|\,.
\end{equation}
Moreover, for any $x,y\in N$ and $\xi\in X$ we have that
\begin{align}
|\xi|_x-|\xi|_y&=
\sup\limits_{t\in\mathbb{S}_{\vartheta_0,0}}\big\|e^{-\beta |t|}(\dd_x S_t)(\xi)\big\|-
\sup\limits_{t\in\mathbb{S}_{\vartheta_0,0}}\big\|e^{-\beta |t|}(\dd_y S_t)(\xi)\big\|\nonumber\\
&\le\sup_{t\in\mathbb{S}_{\vartheta_0,0}}e^{-\beta |t|}\big\|(\dd_x S_t)(\xi)-(\dd_y S_t)(\xi)\big\|\nonumber\\
&\le M\|x-y\| \|\xi\|,\label{eq:Finsler_continuity}
\end{align}
where we used the second estimate in \eqref{eq:dS_t-exponential_growth} and the triangle inequality
for the $L^\infty$-norm on $C_b\big(\mathbb{S}_{\vartheta_0,0},X\big)$.
In particular, the inequality \eqref{eq:Finsler_continuity} implies that the function $F : TN\to\R_{\ge 0}$, $
F(x,\xi):=|\xi|_x$, where $TN=N\times X$ is the tangent bundle of $N$, is continuous. 
Hence, \eqref{eq:Finsler-Euclidean_norms} defines a Finsler metric on $N$.
Now, we define the distance on $N$ in a standard way
\begin{equation}\label{eq:distance}
d(x,y)=d_N(x,y):=\inf_{\gamma}\int_0^1|\dt\gamma(s)|_{\gamma(s)}\,ds
\end{equation}
where the infimum is taken over all 
\begin{equation}\label{eq:gamma_class}
\gamma\in C^1\big([0,1],N\big)\,\,\text{\rm such that}\,\,\gamma(0)=x\,\,\text{\rm and }\,\gamma(1)=y. 
\end{equation}
For any given $t\in\mathbb{S}_{\vartheta_0,0}$, $x,y\in N$, 
and $\gamma$ that satisfies \eqref{eq:gamma_class}, we have
\begin{align}
d\big(S_t(x),S_t(y)\big)&\le\int_0^1\big|\big(S_t(\gamma(s))\big)^\bfdot\big|_{S_t(\gamma(s))}\,ds=
\int_0^1\sup\limits_{\tau\in\mathbb{S}_{\vartheta_0,0}}
\big\|e^{-\beta|\tau|}(\dd_{S_t(\gamma(s))}S_\tau)\big(S_t(\gamma(s))\big)^\bfdot\big\|\,ds\nonumber\\
&=\int_0^1\sup\limits_{\tau\in\mathbb{S}_{\vartheta_0,0}}
\big\|e^{-\beta|\tau|}(\dd_{\gamma(s)}S_{\tau+t})({\dt\gamma}(s))\big\|\,ds\nonumber\\
&\le e^{\beta|t|}\int_0^1\sup\limits_{\tau\in\mathbb{S}_{\vartheta_0,0}}
\big\|e^{-\beta|\tau+t|}(\dd_{\gamma(s)}S_{\tau+t})({\dt\gamma}(s))\big\|\,ds\nonumber\\
&\le e^{\beta|t|}\int_0^1|{\dt\gamma}(s)|_{\gamma(s)}\,ds\label{eq:S-quasi-contraction_preliminary}
\end{align}
where $(\cdot)^\bfdot$ denotes the derivative with respect to $s$.
In the inequalities above we used the fact that $S_t(N)\subseteq N$ and that 
$t+\mathbb{S}_{\vartheta_0,0}\subseteq\mathbb{S}_{\vartheta_0,0}$ for any given $t\in\mathbb{S}_{\vartheta_0,0}$.
It follows from \eqref{eq:Finsler-Euclidean_norms}, the definition \eqref{eq:distance}, 
and \eqref{eq:straight_interval_property}, that
\begin{equation}\label{eq:Finser-Euclidean_distance_preliminary}
\|x-y\|\le d(x,y)\le M\inf_{\gamma}\int_0^1\|\dt\gamma(s)\|\,ds
\end{equation}
for any $x,y\in N$ and $\gamma$ satisfying \eqref{eq:gamma_class}. 
Let $B\subseteq N$ be a non-empty open {\em convex} subset in $N$.
Then, it follows from \eqref{eq:S-quasi-contraction_preliminary}, \eqref{eq:Finser-Euclidean_distance_preliminary},
and Lemma \ref{lem:straight_interval_property} below, that for any $t\in\mathbb{S}_{\vartheta_0,0}$ and 
$x,y\in B$ we have that
\begin{equation}\label{eq:S-quasi-contraction}
d\big(S_t(x),S_t(y)\big)\le e^{\beta|t|}\|x-y\|\quad\text{\rm and}\quad
\|x-y\|\le d(x,y)\le M\|x-y\|\,.
\end{equation}
This proves the first inequality in \eqref{eq:S,G-quasi-contraction}.

Let us now prove the second inequality in \eqref{eq:S,G-quasi-contraction}.
Take $x_0\in N_1\subseteq N$. We then use the continuity of $G : \mathbb{S}_{\vartheta_0,0,T}\times N_1\to X$ 
and Lemma \ref{lem:minimizing_sequence} below to choose open balls 
$U(x_0)\subseteq\widetilde{U}(x_0)\subseteq N_1$ in $X$ centered at $x_0$ and a constant $0<T_0\le T$ such that 
\begin{equation}\label{eq:T_0-relation}
G_t\big(\widetilde{U}(x_0)\big)\subseteq N_1\quad\forall t\in\mathbb{S}_{\vartheta_0,T_0}
\end{equation}
and for any $x,y\in U(x_0)$ there exists a minimizing sequence of $(\gamma_n)_{n\ge 1}$ of 
the action functional $\int_0^1|\dt\gamma(s)|_{\gamma(s)}\,ds$ such that 
\begin{equation}\label{eq:gamma_subclass}
\gamma_n\in C^1\big([0,1],\widetilde{U}(x_0)\big),\,\,\gamma_n(0)=x,\,\text{\rm and}\,\,\gamma_n(1)=y.
\end{equation}
We will complete the proof by modifying the argument in the proof of \cite[Proposition 5.1]{McOwenTopalov6}
(cf. p. 69 in \cite{Marsden}).
For any $x\in N_1$, $\xi\in X$, and $t,s\in\mathbb{S}_{\vartheta_0,\tau_0}$, where $0<\tau_0\le T$ is chosen
so that $G_t(x)\in N_1$ for any $t\in\mathbb{S}_{\vartheta_0,\tau_0}$, we have
\begin{equation}\label{eq:eqlt6}
(\dd_xG_t)(\xi)-(\dd_xG_s)(\xi)=\int_s^t\frac{\dd}{\dd\lambda}(\dd_xG_\lambda)(\xi)\,d\lambda=
\int_s^t(\dd_{G_\lambda(x)}\E)(\dd_xG_\lambda)(\xi)\,d\lambda ,
\end{equation}
where the integration is taken over the straight segment $[s,t]\subseteq\mathbb{S}_{\vartheta_0,0,\tau_0}$ and
where we used that $G : \mathbb{S}_{\vartheta_0,0,\tau_0}\times N_1\to X$ is $C^2$ on conic sets and hence
$\frac{\dd}{\dd\lambda}(\dd_xG_\lambda)(\xi)=(\dd_{G_\lambda(x)}\E)(\dd_xG_\lambda)(\xi)$.
We now take $t\in\mathbb{S}_{\vartheta_0,\tau_0}$ and choose $s\in[0,t]$.
Since the integrand in \eqref{eq:eqlt6} is a continuous function of $\lambda\in[s,t]$ (with values in $X$), 
we obtain from \eqref{eq:E-estimates} and \eqref{eq:Finsler-Euclidean_norms} that
\begin{align}\label{eq:eqlt7}
\Big|\big\3(\dd_xG_t)(\xi)\big\3_{G_s(x)}-\big\3(\dd_xG_s)(\xi)\big\3_{G_s(x)}\Big|
&\le\int_{|s|}^{|t|}\big\3(\dd_{G_{\lambda(\mu)}(x)}\E)(\dd_xG_{\lambda(\mu)})(\xi)\big\3_{G_s(x)}\,
|d\lambda(\mu)|\nonumber\\
&\le ML\int_{|s|}^{|t|}\big\|(\dd_xG_{\lambda(\mu)})(\xi)\big\|\,d\mu ,
\end{align}
where $\lambda(\mu):=\big(\frac{t}{|t|}\big)\mu$ with $\mu\in\big[|s|,|t|\big]$.
On the other side, we obtain from \eqref{eq:Finsler_continuity} that
\begin{align}\label{eq:eqlt8}
\Big|\big\3(\dd_xG_t)(\xi)\big\3_{G_s(x)}-\big\3(\dd_xG_t)(\xi)\big\3_{G_t(x)}\Big|&\le
M\big\|G_s(x)-G_t(x)\big\|\,\big\|(\dd_xG_t)(\xi)\big\|\nonumber\\
&\le ML\big||t|-|s|\big|\,\big\|(\dd_xG_t)(\xi)\big\|
\end{align} 
where we used that 
$\big\|G_s(x)-G_t(x)\big\|=\Big\|\int_{|s|}^{|t|}\E\big(G_{\lambda(\mu)}(x)\big)\,d\lambda(\mu)\Big\|\le L\big||t|-|s|\big|$.
By combining \eqref{eq:eqlt7} and \eqref{eq:eqlt8}, we then conclude that
\begin{equation}\label{eq:eqlt9}
\Big|\big\3(\dd_xG_t)(\xi)\big\3_{G_t(x)}-\big\3(\dd_xG_s)(\xi)\big\3_{G_s(x)}\Big|\le
ML\Big(\int_{|s|}^{|t|}\big\|(\dd_xG_{\lambda(\mu)})(\xi)\big\|\,d\mu+\big||t|-|s|\big|\,\big\|(\dd_xG_t)(\xi)\big\|\Big)
\end{equation}
for any $\xi\in X$, $t\in\mathbb{S}_{\vartheta_0,\tau_0}$, and $s\in[0,t]$. 
Take $t\in\mathbb{S}_{\vartheta_0,\tau_0}$ and let $0=t_1\le t_2\le...\le t_\ell=|t|$ be a partition of 
the interval $\big[0,|t|\big]$. We then obtain from \eqref{eq:eqlt9} that
\begin{align*}
\Big|\big\3(\dd_xG_t)(\xi)\big\3_{G_t(x)}(\xi)-\3\xi\3_x\Big|&\le
\sum_{j=1}^{\ell-1}\Big|\big\3(\dd_xG_{\lambda(t_j)})(\xi)\big\3_{G_{\lambda(t_j)}(x)}-
\big\3(\dd_xG_{\lambda(t_{j+1})}(\xi)\big\3_{G_{\lambda(t_{j+1})}(x)}\Big|\\
&\le ML\Big(\int_0^{|t|}\big\|(\dd_xG_{\lambda(\mu)})(\xi)\big\|\,d\mu+
\sum_{j=1}^{\ell-1}\big\|(\dd_xG_{\lambda(t_{j+1})}(\xi)\big\|\,|t_{j+1}-t_j|
\Big)
\end{align*}
for any $\xi\in X$, and $t\in\mathbb{S}_{\vartheta_0,\tau_0}$.
By passing to the limit as the diameter of the partition $\max\limits_{1\le j\le\ell-1}|t_{j+1}-t_j|\to 0$
we then obtain that
\[
\Big|\big\3(\dd_xG_t)(\xi)\big\3_{G_t(x)}(\xi)-\3\xi\3_x\Big|\le 
2ML\int_0^{|t|}\big\|(\dd_xG_{\lambda(\mu)})(\xi)\big\|\,d\mu
\le 2ML\int_0^{|t|}\big\3(\dd_xG_{\lambda(\mu)})(\xi)\big\3_{G_{\lambda(\mu)}(x)}\,d\mu
\]
where we used \eqref{eq:Finsler-Euclidean_norms} and the fact that the Riemann sums converge since
the integrand is a continuous function on $\big[0,|t|\big]$. This implies that
\[
\big\3(\dd_xG_t)(\xi)\big\3_{G_t(x)}\le
\3\xi\3_x+2ML\int_0^{|t|}\big\3(\dd_xG_{\lambda(\mu)})(\xi)\big\3_{G_{\lambda(\mu)}(x)}\,d\mu
\]
for any $\xi\in X$, and $t\in\mathbb{S}_{\vartheta_0,\tau_0}$.
By Gronwall's inequality,
\begin{equation}\label{eq:eqlt10}
\big\3(\dd_xG_t)(\xi)\big\3_{G_t(x)}\le\3\xi\3_x e^{2ML|t|}
\end{equation}
for any $x\in N_1$, $\xi\in X$, and $t\in\mathbb{S}_{\vartheta_0,\tau_0}$ where $0<\tau_0\le T$ is chosen
so that $G_t(x)\in N_1$ for any $t\in\mathbb{S}_{\vartheta_0,\tau_0}$.
(We apply Gronwall's inequality for
each radial line $t=a\mu$, $\mu\in[0,\tau)$ with fixed direction $a\in\mathbb{S}_{\vartheta_0}$, $|a|=1$.)
Take an open set $V\subseteq U(x_0)$ and choose $0<\tau\le T_0$ such that $G_t(V)\subseteq U(x_0)$ for any 
$t\in\mathbb{S}_{\vartheta_0,\tau}$.
We then obtain that for any $t\in\mathbb{S}_{\vartheta_0,\tau}$, $x,y\in V$, and for any $C^1$-curve 
$\gamma : [0,1]\to\widetilde{U}(x_0)\subseteq N_1$, $\gamma(0)=x$, $\gamma(1)=y$, one has
\begin{equation}\label{eq:eqlt11}
d\big(G_t(x),G_t(y)\big)\le
\int_0^1\big\3\big(G_t\big(\gamma(s)\big)\big)^\bfdot\big\3_{G_t(\gamma(s))}\,ds
\le e^{2ML|t|}\int_0^1\3\dt{\gamma}(s)\3_{\gamma(s)}\,dt
\end{equation}
where we used \eqref{eq:distance} and $G_t\big(\widetilde{U}(x_0)\big)\subseteq N_1\subseteq N$ to conclude 
the inequality appearing on the left side of \eqref{eq:eqlt11}, as well as Gronwall's inequality \eqref{eq:eqlt10} to 
conclude the inequality on the right side.
By taking the minimizing sequence $(\gamma_n)_{n\ge 1}$ of the action functional 
$\int_0^1\3\dt{\gamma}(s)\3_{\gamma(s)}\,dt$ as in \eqref{eq:gamma_subclass} and 
then passing to the limit in \eqref{eq:eqlt11} (with $\gamma$ replaced by $\gamma_n$) as $n\to\infty$,
we conclude that
\[
d\big(G_t(x),G_t(y)\big)\le e^{2ML|t|}d(x,y)
\]
for any $x,y\in V$ and $t\in\mathbb{S}_{\vartheta_0,0,\tau}$.  (The case when $t=0$ is trivial.)
This proves the second inequality in \eqref{eq:S,G-quasi-contraction}.
Since $U(x_0)$ is convex, we obtain from \eqref{eq:S-quasi-contraction} that
the restriction of the distance function \eqref{eq:distance} to $U(x_0)\subseteq N$ 
is equivalent to the standard distance on $U(x_0)$. This completes the proof of the Proposition.

\endproof

The Lemma below is well-known and can be proved in a straightforward way.
\begin{Lem}\label{lem:straight_interval_property}
Let $N\subseteq X$ be an non-empty open set in a Banach space $(X,\|\cdot\|)$. Then,
for any $x,y\in N$ we have that
\begin{equation}\label{eq:straight_interval_property}
\|x-y\|\le\inf\limits_{\gamma}\int_0^\ell\|\dt\gamma(s)\|\,ds,\quad\ell>0,
\end{equation}
where the infimum is taken over all $\gamma\in C^1\big([0,\ell],N\big)$ such that $\gamma(0)=x$ and $\gamma(\ell)=y$.
If the straight segment $[x,y]\subseteq N$ then $\inf\limits_{\gamma}\int_0^\ell\|\dt\gamma(s)\|\,ds=\|x-y\|$.
\end{Lem}

The Lemma below is used in the proof of Lemma \ref{lem:equivalent_distance}.
It shows that a minimizing sequence of the action $\int_0^1|\dt\gamma(s)|_{\gamma(s)}\,ds$
can be chosen ``not too far away'' (with respect to the norm in $X$) from the points $x,y\in N$.

\begin{Lem}\label{lem:minimizing_sequence}
Let $B_\rho^{\|\cdot\|}(x_0)\subseteq N$ be an open ball of radius $\rho>0$ centered
at $x_0\in N$ where $N$ is the open set in $X$ appearing in Lemma \ref{lem:equivalent_distance}. 
Then, there exists $r=r(\rho)>0$ such that for any $x,y\in B_r^{\|\cdot\|}(x_0)$
there exists a minimizing sequence of curves $(\gamma_n)_{n\ge 0}$,
$\gamma_n\in C^1\big([0,1],B_\rho^{\|\cdot\|}(x_0)\big)$, $\gamma_n(0)=x$, $\gamma_n(1)=y$, $n\ge 1$,
such that
\[
d(x,y)=\lim_{n\to\infty}\int_0^1|\dt\gamma_n(s)|_{\gamma_n(s)}\,ds
\]
where $d(x,y)$ is the distance \eqref{eq:distance} on $N$.
\end{Lem}

\begin{proof}[Proof of Lemma \ref{lem:minimizing_sequence}]
Choose $r:=\rho/4M$ where $M\ge 1$ is the constant in \eqref{eq:S-quasi-contraction}. 
Then, for any $x,y\in B_r^{\|\cdot\|}(x_0)$ we obtain from
the triangle inequality and the second estimate in \eqref{eq:S-quasi-contraction} that
\begin{equation}\label{eq:d<rho}
d(x,y)\le d(x,x_0)+d(y,x_0)\le M\big(\|x-x_0\|+\|y-x_0\|\big)<\rho/2\,.
\end{equation}
Now, take an arbitrary minimizing sequence $(\gamma_n)_{n\ge 1}$ of the action functional
$\int_0^1|\dt\gamma(s)|_{\gamma(s)}\,ds$ that satisfies \eqref{eq:gamma_class}.
Since $d(x,y)=\lim\limits_{n\to\infty}\int_0^1|\dt\gamma_n(s)|_{\gamma_n(s)}\,ds$, it
follows from \eqref{eq:d<rho} that we can choose $n_*\ge 1$ such that 
\[
\int_0^1|\dt\gamma_n(s)|_{\gamma_n(s)}\,ds<\rho/2
\]
for any $n\ge n_*$.
It now follows from \eqref{eq:Finsler-Euclidean_norms} and \eqref{eq:straight_interval_property} that
\[
\|x-\gamma_n(t)\|\le\int_0^t\|\dt\gamma_n(s)\|\,ds\le\int_0^1|\dt\gamma_n(s)|_{\gamma_n(s)}\,ds<\rho/2
\]
for any $t\in[0,1]$ and for any $n\ge n_*$. This implies that $\gamma_n(t)\in B_\rho^{\|\cdot\|}(x_0)$ for 
any $t\in[0,1]$ and for any $n\ge n_*$. By rearranging the indexes of the sequence $(\gamma_n)_{n\ge n_*}$,
we then complete the proof of the Lemma.
\end{proof}

The following Lemma is crucial for the proof of Proposition \ref{prop:product_formula}. 


\begin{Lem}\label{lem:stability}
Let $(Z,\|\cdot\|_Z)$ and $(X,\|\cdot\|_X)$ be complex Banach spaces 
such that the inclusion $Z\subseteq X$ is bounded, $Z$ is dense in $X$, and let $U$ be a non-empty open set in $X$. 
Assume that the maps 
\begin{equation}\label{eq:G,S}
G : \mathbb{S}_{\vartheta_0,0,T}\times U\to X
\quad\text{\rm and}\quad 
S : \mathbb{S}_{\vartheta_0,0,T}\times U\to X
\end{equation}
are continuous, satisfy the local semigroup property in $U$, and they are locally equivalent to a quasi-contraction in $U$ with 
the same distance function. Assume in addition that the maps
\begin{equation}\label{eq:G,S-diagonal}
G : \mathbb{S}_{\vartheta_0,0,T}\times U_Z\to X
\quad\text{\rm and}\quad 
S : \mathbb{S}_{\vartheta_0,0,T}\times U_Z\to X
\end{equation}
are $C^1$ on conic sets. (This is summarized in diagram \eqref{eq:diagram1}.)
\begin{equation}\label{eq:diagram1}
\begin{tikzcd}
U\arrow[r, "S_t{,}\,G_t"]&X\\
U_Z\arrow[u, hook]\arrow[ur, swap, "S_t{,}\,G_t"]&         
\end{tikzcd}
\end{equation}
Then, for any $x_0\in U$ and for any open neighborhood $\widetilde{V}(x_0)\subseteq U$ of $x_0$ in $X$ there exist 
an open neighborhood $V(x_0)\subseteq\widetilde{V}(x_0)$ of $x_0$ in $X$ and $0<\tau=\tau(x_0)\le T$ 
(generally, depending on the choice of $x_0$ and $\widetilde{V}(x_0)$) such that for any $x\in V(x_0)$ and 
$t\in\mathbb{S}_{\vartheta_0,0,\tau}$ we have that
\begin{equation}\label{eq:stability}
S_{t/n}^{\alpha_1}\circ G_{t/n}^{\mu_1}\circ\cdots\circ S_{t/n}^{\alpha_\ell}\circ 
G_{t/n}^{\mu_\ell}(x)\in\widetilde{V}(x_0)
\end{equation}
for any $n\ge 1$ and for any multi-indexes $\alpha,\mu\in\Z_{\ge 0}^\ell$ with $1\le\ell\le 2n$,
$|\alpha|\le n$, $|\mu|\le n$, and where we set $S_{t/n}^0(x)=x$ and $G_{t/n}^0(x)=x$ for $x\in U$.
\end{Lem}

\begin{Rem}\label{rem:stability}
In particular, Lemma \ref{lem:stability} implies that for any $x\in V(x_0)$ and $t\in\mathbb{S}_{\vartheta_0,0,\tau}$ 
we have that $(S_{t/n}\circ G_{t/n})^n(x)\in\widetilde{V}(x_0)$ for any $n\ge 1$.
\end{Rem}

\begin{proof}[Proof of Lemma \ref{lem:stability}] 
Take $x_0\in U$. Since the maps in \eqref{eq:G,S} satisfy condition $1)$ in $(b)$, we can assume without loss of generality that
there is a distance function $d : U\times U\to\R_{\ge 0}$ equivalent to the standard distance on $U$ so that
\eqref{eq:equivalent_distances} holds and $\widetilde{V}(x_0)\equiv U$. 
By the continuity of the maps in \eqref{eq:G,S}, there exist an open neighborhood $V(x_0)\subseteq U$ of $x_0$ in $X$ and
a constant $0<\tau=\tau(x_0)\le T$ such that 
\begin{equation}\label{eq:G,S_in_U}
S_t(x), G_t(x)\in U
\end{equation}
for any $x\in V(x_0)$ and $t\in\mathbb{S}_{\vartheta_0,0,\tau}$.
It also follows from the local semigroup property that it is sufficient to prove that
for any $x\in V(x_0)$ and $t\in\mathbb{S}_{\vartheta_0,0,\tau}$ 
(with $V(x_0)$ and $\tau(x_0)>0$ taken smaller if necessary)
we have that
\begin{equation}\label{eq:stability'}
S_{t/n}^{\alpha_1}\circ G_{t/n}^{\mu_1}\circ\cdots\circ S_{t/n}^{\alpha_\ell}\circ 
G_{t/n}^{\mu_\ell}(x)\in U
\end{equation}
for any $n\ge 1$, $1\le\ell\le 2n$, $|\alpha|\le n$, $|\mu|\le n$ and 
\begin{equation}\label{eq:reduced_condition}
\alpha_j,\mu_j\in\{0,1\},\quad\alpha_j\ne\mu_j,\quad 1\le j\le\ell\,.
\end{equation}
Let us prove that. For given $z_0\in V(x_0)\cap Z$ and $\lambda_1,\lambda_2\in\{0,1\}$, $\lambda_1\ne\lambda_2$, 
consider the map $\mathbb{S}_{\vartheta_0,0,\tau}\to U$, 
$t\mapsto z_{\lambda_1,\lambda_2}(t):=S_{\lambda_1 t}\circ G_{\lambda_2 t}(z_0)$.
Note that $z_{1,0}(t)=S_t(z_0)$ and $z_{0,1}(t)=G_t(z_0)$.
This and the assumptions on the maps \eqref{eq:G,S-diagonal} imply that the map 
$\mathbb{S}_{\vartheta_0,0,\tau}\to U$, $t\mapsto z_{\lambda_1,\lambda_2}(t)$,
is $C^1$, and hence, for $t\in\mathbb{S}_{\vartheta_0,0,\tau}$,
\[
\|z_{\lambda_1,\lambda_2}(t)-z_0\|\le\Big\|t\int_0^1\frac{\dd z_{\lambda_1,\lambda_2}}{\dd t}(s t)\,ds\Big\|\le K|t|,
\]
where $\|\cdot\|$ is the norm in $X$, and where we used the continuity of the map
\[
\mathbb{S}_{\vartheta_0,0,\tau}\to X,\quad t\mapsto\frac{\dd z_{\lambda_1,\lambda_2}}{\dd t}(t),
\] 
at $t=0$ (see Definition \ref{def:C^k_along_curves_complex} and Definition \ref{def:C^k_complex}
in Appendix \ref{sec:C^k_along_curves}) to chose $K>0$ and $\tau>0$ (taken smaller if necessary) so that 
\[
\sup\limits_{t\in\mathbb{S}_{\vartheta_0,0,\tau}}\Big\|\frac{\dd z_{\lambda_1,\lambda_2}}{\dd t}(t)\Big\|\le K<\infty
\]
for any $\lambda_1,\lambda_2\in\{0,1\}$, $\lambda_1\ne\lambda_2$.
This and \eqref{eq:equivalent_distances} then imply that for any $t\in\mathbb{S}_{\vartheta_0,0,\tau}$ and 
for any $\lambda_1,\lambda_2\in\{0,1\}$, $\lambda_1\ne\lambda_2$, we have that
\begin{equation}\label{eq:est1}
d\big(z_0, S_{\lambda_1 t}\circ G_{\lambda_2 t}(z_0)\big)\le K_1|t|
\end{equation}
where $K_1:=K/L_1$.
Let $\delta>0$ be chosen so small that the ball $B^d_\delta(z_0):=\{x\in X\,|\,d(x,z_0)<\delta\}$ is 
contained in $V(x_0)$. By taking $\tau>0$ smaller if necessary, we ensure that
\begin{equation}\label{eq:est2}
2e^{2\beta\tau}K_1\tau<\delta
\end{equation}
where $\beta>0$ is the constant appearing in \eqref{eq:quasi-contraction}.
We now take $n\ge 1$, $\alpha,\mu\in\Z_{\ge 0}^\ell$, $1\le\ell\le 2n$, $|\alpha|\le n$, $|\mu|\le n$, 
such that the condition \eqref{eq:reduced_condition} holds.
We will first prove that for any $n\ge 1$, $1\le j\le\ell$, $1\le\ell\le 2n$, and $t\in\mathbb{S}_{\vartheta_0,0,\tau}$, 
we have that
\begin{equation}\label{eq:est3}
d\big(z_0,\Pi_j(z_0)\big)<\delta
\end{equation}
where
\[
\Pi_j(z_0):=(S_{t/n}^{\alpha_1}\circ G_{t/n}^{\mu_1})\circ\cdots
\circ(S_{t/n}^{\alpha_j}\circ G_{t/n}^{\mu_j})(z_0)\in V(x_0)\,.
\] 
To this end, note that it follows from \eqref{eq:est1}
(applied with $\lambda_1=\alpha_j$, $\lambda_2=\mu_j$, and $t$ replaced by $t/n$) 
and \eqref{eq:est2} that for any $t\in\mathbb{S}_{\vartheta_0,0,\tau}$,
\begin{equation}\label{eq:step1}
d\big(z_0,(S_{t/n}^{\alpha_j}\circ G_{t/n}^{\mu_j})(z_0)\big)\le K_1 |t|/n<\delta,
\quad 1\le j\le\ell,
\end{equation}
and hence $(S_{t/n}^{\alpha_j}\circ G_{t/n}^{\mu_j})(z_0)\in V(x_0)$, $1\le j\le\ell$. 
This proves the statement when $\ell=1$ and when $j=1$ and $2\le\ell\le 2n$.
Let us now assume that $2\le\ell\le 2n$.
Then, we obtain from \eqref{eq:G,S_in_U} that $\Pi_2(z_0)\in U$ for any $t\in\mathbb{S}_{\vartheta_0,0,\tau}$, 
and hence,
\begin{align*}
d\big(z_0,\Pi_2(z_0)\big)&\le
d\big(z_0,(S_{t/n}^{\alpha_1}\circ G_{t/n}^{\mu_1})(z_0)\big)+
d\big((S_{t/n}^{\alpha_1}\circ G_{t/n}^{\mu_1})(z_0),\Pi_2(z_0)\big)\\
&\le d\big(z_0,(S_{t/n}^{\alpha_1}\circ G_{t/n}^{\mu_1})(z_0)\big)
+e^{\beta(\alpha_1+\mu_1)|t|/n}d\big(z_0,(S_{t/n}^{\alpha_2}\circ G_{t/n}^{\mu_2})(z_0)\big)\\
&\le 2e^{\beta\tau}K_1 |t|/n\le 2e^{2\beta\tau}K_1\tau<\delta,
\end{align*}
where we used the property \eqref{eq:quasi-contraction}, the first inequality in \eqref{eq:step1}, and \eqref{eq:est2}.
In particular, we see that $\Pi_2(z_0)\in V(x_0)$ for $t\in\mathbb{S}_{\vartheta_0,0,\tau}$, 
and therefore $\Pi_3(z_0)\in U$.
By continuing this process inductively in $2\le j\le\ell$ (with $\ell\le 2n$ fixed), we obtain that for any 
$t\in\mathbb{S}_{\vartheta_0,0,\tau}$,
\begin{align*}
d\big(z_0,\Pi_j(z_0)\big)&\le
\sum_{1\le k\le j}d\big(\Pi_{k-1}(z_0),\Pi_k(z_0)\big)\\
&\le\sum_{1\le k\le j}e^{\beta(k-1)|t|/n}d\big(z_0,(S_{t/n}^{\alpha_k}\circ G_{t/n}^{\nu_k})(z_0)\big)\\
&\le j e^{\beta(j-1)|t|/n} K_1 |t|/n\le 2e^{2\beta\tau}K_1\tau<\delta,
\end{align*}
where we used \eqref{eq:quasi-contraction}, \eqref{eq:step1}, \eqref{eq:est2}, the fact that $\alpha_q,\mu_q\in\{0,1\}$, 
$\alpha_q\ne\mu_q$, $1\le q\le\ell$, and the inequality $j\le\ell\le 2n$. Hence, $\Pi_j(z_0)\in V(X_0)$ for 
$t\in\mathbb{S}_{\vartheta_0,0,\tau}$. This proves the preliminary statement \eqref{eq:est3}.
We now take $\epsilon>0$ so that $B^d_{3\epsilon}(x_0)\subseteq V(x_0)$ and choose
$z_0\in B^d_\kappa(x_0)\cap Z$ where $0<\kappa\le\epsilon$ will be chosen below. 
By the discussion above, there exists $\tau>0$ (depending on the choice of $z_0$)
so that \eqref{eq:est3} (with $\delta=\kappa$) holds for any $n\ge 1$, $1\le\ell\le 2n$, 
$|\alpha|\le n$, $|\beta|\le n$, $1\le j\le\ell$, and $t\in\mathbb{S}_{\vartheta_0,0,\tau}$. 
Let us now choose $\kappa>0$ so that
\begin{equation}\label{eq:est4}
2 e^{2\beta\tau}\kappa<\epsilon\,.
\end{equation}
Then, for any $x\in B^d_\kappa(x_0)$ and for any $n\ge 1$ we argue inductively in $1\le j\le\ell$
to conclude that 
\[
\Pi_j(x):=(S_{t/n}^{\alpha_1}\circ G_{t/n}^{\mu_1})\circ\cdots
\circ(S_{t/n}^{\alpha_j}\circ G_{t/n}^{\mu_j})(x)\in U
\]
and that for any $t\in\mathbb{S}_{\vartheta_0,0,\tau}$,
\begin{align*}
d\big(x_0,\Pi_j(x)\big)&\le
d(x_0,z_0)+d\big(z_0,\Pi_j(z_0)\big)+d\big(\Pi_j(z_0),\Pi_j(x)\big)\\
&\le 2\kappa+2e^{\beta j |t|/n}\kappa\le 2\epsilon+2e^{2\beta\tau}\kappa<3\epsilon,
\end{align*}
where we used \eqref{eq:quasi-contraction}, \eqref{eq:est3} (with $\delta=\kappa$), and \eqref{eq:est4}.
In particular, $\Pi_j(x)\in V(x_0)$, and hence $\Pi_{j+1}\in U$.
Finally, we complete the proof by redefining $V(x_0):=B^d_\kappa(x_0)$.
\end{proof}

\begin{Rem}\label{rem:nu}
The proof of Lemma \ref{lem:stability} shows that for any $\widetilde{V}(x_0)\subseteq U$ there exist 
an open neighborhood $V=V(x_0)\subseteq\widetilde{V}(x_0)$ of $x_0$ in $X$ and 
$0<\tau\le T$ as above such that \eqref{eq:stability} holds with $S_t$ replaced by $S^{(\nu)}_t=S_{\nu t}$ 
uniformly in $0\le\nu\le 1$. This follows since the second estimate in \eqref{eq:quasi-contraction}
with $S_t$ replaced by $S^{(\nu)}_t$ and the estimate 
$\sup\limits_{t\in\mathbb{S}_{\vartheta_0,0,\tau}}\big\|\frac{\dd z^{(\nu)}}{\dd t}(t)\big\|\le K$ where 
$z^{(\nu)}(t):=S^{(\nu)}_{\lambda_1 t}\circ G_{\lambda_2 t}(z_0)$, $z_0\in U_Z$, with $\lambda_1,\lambda_2\in\{0,1\}$,
$\lambda_1\ne\lambda_2$, hold uniformly in $0\le\nu\le 1$.
\end{Rem}

\begin{proof}[Proof of Proposition \ref{prop:product_formula}]
We first note that for any given value of $0\le\nu\le 1$ the map $S^{(\nu)} : \mathbb{S}_{\vartheta_0,0,T}\times U\to X$
where $S_t^{(\nu)}(x)=S_{\nu t}(x)$, $x\in U$, is well-defined and satisfies all the assumptions on
the map $S : \mathbb{S}_{\vartheta_0,0,T}\times U\to X$ in Proposition \ref{prop:product_formula}.
It follows from condition $(a)$, condition $1)$ in $(b)$, and Lemma \ref{lem:stability} (and Remark \ref{rem:nu}), 
that for any $x_0\in U$ we can choose an open ball $B\subseteq U$ in $X$ centered at $x_0$ whose
closure in $X$ is contained in $U$, an open neighborhood $O(x_0)\subseteq B$ of $x_0$ in $X$, 
and constants $0<\tau=\tau(x_0)<T$ and $\beta>0$, such that 
there exists a distance function $d : B\times B\to\R_{\ge 0}$ equivalent to the standard distance on $B$
and such that:
\begin{itemize}
\item[$\bullet$] For any $x\in O(x_0)$, $t\in\mathbb{S}_{\vartheta_0,0,\tau}$, and $0\le\nu\le 1$, we have that
\begin{equation}\label{eq:stability_nu}
P_{t,n,\nu}^{\alpha,\mu}(x):=
\big(S_{t/n}^{(\nu)}\big)^{\alpha_1}\circ G_{t/n}^{\mu_1}\circ\cdots\circ
\big(S_{t/n}^{(\nu)}\big)^{\alpha_\ell}\circ G_{t/n}^{\mu_\ell}(x)\in B
\end{equation}
for any $n\ge 1$ and for any multi-indexes $\alpha,\mu\in\Z_{\ge 0}^\ell$ with $1\le\ell\le 2n$,
$|\alpha|\le n$, $|\mu|\le n$, and where we set $S_{t/n}^0(x)=x$ and $G_{t/n}^0(x)=x$ for $x\in O(x_0)$.

\item[$\bullet$] For any $x,y\in O(x_0)$ and $t\in\mathbb{S}_{\vartheta_0,0,\tau}$ we have that 
\begin{equation}\label{eq:quasi-contraction'}
d\big(G_t(x),G_t(y)\big)\le e^{\beta |t|}d(x,y)\quad\text{\rm and}\quad
d\big(S_t^{(\nu)}(x),S_t^{(\nu)}(y)\big)\le e^{\beta |t|}d(x,y)\,.
\end{equation}
\end{itemize}
An immediate consequence from \eqref{eq:stability_nu} and \eqref{eq:quasi-contraction'} is that
for any $x,y\in O(x_0)$ and $t\in\mathbb{S}_{\vartheta_0,0,\tau}$,
\begin{equation}\label{eq:pre-Lipschitz}
d\big(P_{t,n,\nu}^{\alpha,\mu}(x),P_{t,n,\nu}^{\alpha,\mu}(y)\big)\le 
e^{\beta(\alpha|+|\mu|)|t|/n}d(x,y)\le e^{2\beta |t|}d(x,y)
\end{equation}
for any $n\ge 1$ and for any $\alpha,\mu\in\Z_{\ge 0}^\ell$ with $1\le\ell\le 2n$,
$|\alpha|\le n$, and $|\mu|\le n$.

We will also need the following technical Lemma. Note that this is the only place in 
the proof of Proposition \ref{prop:product_formula} where we use the $C^2$ regularity at $t=0$
(see condition $(a)$).

\begin{Lem}\label{lem:commutator_estimate}
For any $z_0\in O(x_0)\cap Z$ there exist an open neighborhood $V_Z(z_0)\subseteq O(x_0)\cap Z$ of $z_0$ in $Z$,
a constant $0<\tau_1=\tau_1(z_0)\le\tau$, and a constant $K>0$, such that 
\begin{equation}\label{eq:commutator_estimate}
d\big(S_t^{(\nu)}\circ G_t(z),G_t\circ S_t^{(\nu)}(z)\big)\le K|t|^2
\end{equation}
for any $z\in V_Z(z_0)$, $t\in\mathbb{S}_{\vartheta_0,0,\tau_1}$, and $0\le\nu\le 1$.
\end{Lem}

\begin{proof}[Proof of Lemma \ref{lem:commutator_estimate}] 
It follows from condition $(a)$ in Proposition \ref{prop:product_formula} and 
Corollary \ref{coro:C^k_along_curves} in Appendix \ref{sec:C^k_along_curves} that
for any $z\in O(x_0)\cap Z$ and $0\le\nu\le 1$ the map
\[
\mathbb{S}_{\vartheta_0,0,\tau}\to X,
\quad t\mapsto R^{(\nu)}(t):=S_t^{(\nu)}\circ G_t(z)-G_t\circ S_t^{(\nu)}(z)
\]
is well-defined and $C^2$.
By Taylor's formula with reminder in integral form, we obtain that
\[
R^{(\nu)}(t)=R^{(\nu)}(0)+\frac{\dd R^{(\nu)}}{\dd t}(0)\,t+
\Big(\int_0^1(1-s)\frac{\dd^2 R^{(\nu)}}{\dd t^2}(st)\,ds\Big)\,t^2
\]
for any $t\in\mathbb{S}_{\vartheta_0,0,\tau}$ and $0\le\nu\le 1$.
Since by the product rule $R^{(\nu)}(0)=\frac{\dd R^{(\nu)}}{\dd t}(0)=0$ we then
conclude from the weak continuity of the directional partial derivatives of order $\le 2$ 
of the maps in \eqref{eq:G} and the map $S : \mathbb{S}_{\vartheta_0,0,T}\times U_Z\to X$ 
(see Definition \ref{def:C^k_along_curves_complex}) that for any $z_0\in O(x_0)\cap Z$ there exist
an open neighborhood $V_Z(z_0)\subseteq O(x_0)\cap Z$ of $z_0$ in $Z$,
a constant $0<\tau_1=\tau_1(z_0)\le\tau$, and a constant $K_1>0$ such that
\begin{equation}\label{eq:commutator_estimate'}
\big\|S_t^{(\nu)}\circ G_t(z)-G_t\circ S_t^{(\nu)}(z)\big\|\le K_1|t|^2
\end{equation}
for any $z\in V_Z(z_0)$, $t\in\mathbb{S}_{\vartheta_0,0,\tau_1}$, and $0\le\nu\le 1$.
Recall from \eqref{eq:stability_nu} that
\[
S_t^{(\nu)}\circ G_t(z),\,\,G_t\circ S_t^{(\nu)}(z)\in B
\]
for any $z\in V_Z(z_0)$, $t\in\mathbb{S}_{\vartheta_0,0,\tau_1}$, and $0\le\nu\le 1$.
Since by construction the distance function $d$ is equivalent to the standard distance on $B$, 
we then conclude from \eqref{eq:commutator_estimate'} and \eqref{eq:equivalent_distances} that
\eqref{eq:commutator_estimate} holds with $K=K_1/L_1$ for any $z\in V_Z(z_0)$, 
$t\in\mathbb{S}_{\vartheta_0,0,\tau_1}$, and $0\le\nu\le 1$.
\end{proof}

\begin{Rem}
It follows easily from the compactness of the closure of $\mathbb{S}_{\vartheta_0',\tau}$ in $\C$ that 
Lemma \ref{lem:commutator_estimate} holds with $\tau_1=\tau$ and $\vartheta_0$ replaced by
$\vartheta_0'$ such that $0<\vartheta_0'<\vartheta_0$.
\end{Rem}

\noindent The following Lemma is a holomorphic version of \cite[Lemma 5.3]{McOwenTopalov6} (cf. \cite[Lemma 3]{EM}).
Note that this is the only place in the proof of Proposition \ref{prop:product_formula} where  we use 
the $C^1$ regularity of the map $S : \mathbb{S}_{\vartheta_0,0,T}\times U_D\to Z$ at $t=0$ and 
the quasi-contraction condition $2)$ in $(b)$).

\begin{Lem}\label{lem:fundamental}
Take $z_0\in O(x_0)\cap Z$ and choose an open neighborhood $V_Z(z_0)\subseteq O(x_0)\cap Z$ and constants
$\tau_1=\tau_1(z_0)>0$ and $K>0$ as in Lemma \ref{lem:commutator_estimate}.
Then, there exist an open neighborhood $W_Z(z_0)\subseteq V_Z(z_0)$ of $z_0$ in $Z$
and $0<\tau_2=\tau_2(z_0)\le\tau_1$ such that for any  $z\in W_Z(z_0)$, $t\in\mathbb{S}_{\vartheta_0,0,\tau_2}$,
$0\le\nu\le 1$, and  $n\ge 1$, we have that
\begin{equation}\label{eq:stability_Z}
\big(S_{t/n}^{(\nu)}\circ G_{t/n}\big)^n(z)\in V_Z(z_0)
\end{equation}
and 
\begin{equation}\label{eq:fundamental}
d\big(\big(S_{t/n}^{(\nu)}\circ G_{t/n}\big)^n(z),\big(S_{t/m}^{(\nu)}\circ G_{t/m}\big)^m(z)\big)\le 
2 K e^{2\beta|t|}|t|^2/n
\end{equation}
for any $m\ge n$.
\end{Lem}

\begin{proof}[Proof of Lemma \ref{lem:fundamental}]
Take $z_0\in O(x_0)\cap Z$. We then use condition $(a)$, condition $2)$ in $(b)$, and  Lemma \ref{rem:stability}, 
to find an open neighborhood $W_Z(z_0)\subseteq V_Z(z_0)$ of $z_0$ in $Z$ and $0<\tau_2=\tau_2(z_0)\le\tau_1$ such that 
for any $z\in W_Z(z_0)$, $t\in\mathbb{S}_{\vartheta_0,0,\tau_2}$, $0\le\nu\le 1$, and $n\ge 1$, we have that
\begin{equation}\label{eq:stability_Z_general}
\big(S_{t/n}^{(\nu)}\big)^{\alpha_1}\circ G_{t/n}^{\mu_1}\circ\cdots\circ
\big(S_{t/n}^{(\nu)}\big)^{\alpha_\ell}\circ G_{t/n}^{\mu_\ell}(z)\in V_Z(z_0)
\end{equation}
for any multi-indexes $\alpha,\mu\in\Z_{\ge 0}^\ell$ with $1\le\ell\le 2n$, $|\alpha|\le n$, $|\mu|\le n$.
In this way, we can apply the inequality \eqref{eq:commutator_estimate} with $z$ replaced by any element
of the form \eqref{eq:stability_Z_general}.
Let us now prove \eqref{eq:fundamental}. We will first assume that $m=k n$ for 
some integer $k\ge 1$. Then, for any $z\in W_Z(z_0)$,
\[
(S_{t/m}^{(\nu)}\circ G_{t/m})^m(z)=\big[(S_{t/(kn)}^{(\nu)}\circ G_{t/(kn)})^k\big]^n(z)
\]
and, by the semigroup property,
\[
(S_{t/n}^{(\nu)}\circ G_{t/n})^n(z)=\big[(S_{t/(kn)}^{(\nu)})^k\circ G_{t/(kn)}^k\big]^n(z)\,.
\]
Consider for simplicity the case when $k=3$. Then, the term 
\[
(S_{t/(3n)}^{(\nu)}\circ G_{t/(3n)})^3(z)=
(S_{t/(3n)}^{(\nu)}\circ G_{t/(3n)})\circ(S_{t/(3n)}^{(\nu)}\circ G_{t/(3n)})\circ(S_{t/(3n)}^{(\nu)}\circ G_{t/(3n)})(z)
\]
can be transformed into $(S_{t/(3n)}^{(\nu)})^3\circ G_{t/(3n)}^3(z)$ by the sequence of terms
\begin{align*}
&A_1:=S_{t/(3n)}^{(\nu)}\circ G_{t/(3n)}\circ S_{t/(3n)}^{(\nu)}\circ(G_{t/(3n)}\circ S_{t/(3n)}^{(\nu)})\circ G_{t/(3n)}(z)\\
&A_2:=S_{t/(3n)}^{(\nu)}\circ(G_{t/(3n)}\circ S_{t/(3n)}^{(\nu)})\circ S_{t/(3n)}^{(\nu)}\circ G_{t/(3n)}\circ G_{t/(3n)}(z)\\
&A_3:=S_{t/(3n)}^{(\nu)}\circ S_{t/(3n)}^{(\nu)}\circ(G_{t/(3n)}\circ S_{t/(3n)}^{(\nu)})\circ G_{t/(3n)}\circ G_{t/(3n)}(z)\\
&A_4:=S_{t/(3n)}^{(\nu)}\circ S_{t/(3n)}^{(\nu)}\circ S_{t/(3n)}^{(\nu)}\circ G_{t/(3n)}\circ G_{t/(3n)}\circ G_{t/(3n)}(z)
\end{align*}
with the property that each term in the sequence is obtained from the previous one by a binary commutation of 
the type described above. In the case of arbitrary $k\ge 1$, we obtain in the same way a sequence of terms 
\[
A_1,...,A_{\ell_k}
\] 
such that each term is obtained from the previous one by a binary commutation,
\[
A_1:=(S_{t/(kn)}^{(\nu)}\circ G_{t/(kn)})^k(z),\quad A_{\ell_k}:=(S_{t/(kn)}^{(\nu)})^k\circ G_{t/(kn)}^k(z),
\]
and $\ell_k-1\le (k-1)+(k-2)+...+1=k(k-1)/2\le k^2$. Then, by the triangle inequality,
\begin{align}
&d\big((S_{t/n}^{(\nu)}\circ G_{t/n})^n(z),(S_{t/m}^{(\nu)}\circ G_{t/m})^m(z)\big)=\nonumber\\
&=d\Big(\big[(S_{t/(kn)}^{(\nu)})^k\circ G_{t/(kn)}^k\big]^n(z),
\big[(S_{t/(kn)}^{(\nu)}\circ G_{t/(kn)})^k\big]^n(z)\Big)\nonumber\\
&\le\sum_{j=0}^{n-1}\sum_{l=1}^{\ell_k-1}d\big(B_{j,l}(z),B_{j,l+1}(z)\big)\label{eq:est5}
\end{align}
where
\begin{equation}\label{eq:B_{j,l}}
B_{j,l}(z):=
\big[(S_{t/(kn)}^{(\nu)}\circ G_{t/(kn)})^k\big]^{n-j-1}\circ A_l\circ\big[(S_{t/(kn)}^{(\nu)})^k\circ G_{t/(kn)}^k\big]^j(z)\,.
\end{equation}
Note that the terms appearing in \eqref{eq:B_{j,l}} are of the form \eqref{eq:stability_Z_general},
and hence belong to $V_Z(z_0)\subseteq U_Z\subseteq U$ where the distance $d$ is defined.
Moreover, since the terms in \eqref{eq:stability_Z_general} belong to $V_Z(z_0)$,
the estimate \eqref{eq:commutator_estimate} can be applied for $z$ replaced by any of the terms appearing in 
\eqref{eq:stability_Z_general}.
In particular, we obtain from \eqref{eq:quasi-contraction}, \eqref{eq:commutator_estimate}, \eqref{eq:B_{j,l}}, 
and the definition of the terms $A_1,...,A_{\ell_k}$, that for any $t\in\mathbb{S}_{\vartheta_0,0,\tau_2}$,
\begin{align}\label{eq:d(B,B)}
d\big(B_{j,l}(z),B_{j,l+1}(z)\big)&\le e^{2\beta |t|(n-1)k/(kn)}e^{2\beta |t|(k-1)/(kn)} K(|t|/{nk})^2\nonumber\\
&\le K e^{2\beta |t|}(|t|/{nk})^2\,.
\end{align}
This, together with \eqref{eq:est5}, then implies that for $m=kn$,
\begin{equation}\label{eq:crucial}
d\big((S_{t/m}^{(\nu)}\circ G_{t/m})^m(z),(S_{t/n}^{(\nu)}\circ G_{t/n})^n(z)\big)\le (n k^2) K e^{2\beta |t|}(|t|/{nk})^2
\le K e^{2\beta |t|}|t|^2/n
\end{equation}
where we also used that $\ell_k-1\le k^2$. For $1\le n\le m$ arbitrary we then have
\begin{align*}
&d\big((S_{t/n}^{(\nu)}\circ G_{t/n})^n(z),(S_{t/m}^{(\nu)}\circ G_{t/m})^m(z)\big)
\le d\big((S_{t/n}^{(\nu)}\circ G_{t/n})^n(z),(S_{t/(nm)}^{(\nu)}\circ G_{t/(nm)})^{nm}(z)\big)+\\
&+d\big((S_{t/(nm)}^{(\nu)}\circ G_{t/(nm)})^{nm}(z),(S_{t/m}^{(\nu)}\circ G_{t/m})^m(z)\big)
\le K e^{2\beta |t|}|t|^2/n+K e^{2\beta |t|}|t|^2/m\\
&\le 2K e^{2\beta |t|}|t|^2/n\,.
\end{align*}
Hence, the inequality \eqref{eq:fundamental} holds for any $1\le n\le m$,
$z\in W(z_0)$, and $t\in\mathbb{S}_{\vartheta_0,0,\tau_2}$.
\end{proof}

It follows from the inequality \eqref{eq:fundamental} in Lemma \ref{lem:fundamental} that for any given $z\in O(x_0)\cap Z$
there exists $0<\tau_2=\tau_2(z)\le\tau$ such that for any $t\in\mathbb{S}_{\vartheta_0,0,\tau_2}$ and $0\le\nu\le 1$ 
there exists a limit in $X$,
\begin{equation}\label{eq:H-definition}
H_t^{(\nu)}(z):=\lim_{n\to\infty}\big(S_{t/n}^{(\nu)}\circ G_{t/n}\big)^n(z)\in X,
\end{equation}
which is uniform in the parameters $t$ and $\nu$ (and locally uniform in $z\in O(x_0)\cap Z$).
In particular, for any $z\in O(x_0)\cap Z$ and $0\le\nu\le 1$ we can consider the map
\begin{equation}\label{eq:H-curve}
\mathbb{S}_{\vartheta_0,\tau_2}\to X,\quad t\mapsto H_t^{(\nu)}(z),
\end{equation}
where $\tau_2>0$ depends on the choice of $z\in O(x_0)\cap Z$ and $\tau_2\le\tau$.
Recall from \eqref{eq:stability_nu} and the analyticity of the maps \eqref{eq:S,G-restricted}
that for any given values of $z\in O(x_0)\cap Z$ and $0\le\nu\le 1$ the maps
\[
\mathbb{S}_{\vartheta_0,\tau}\to X,\quad t\mapsto\big(S_{t/n}^{(\nu)}\circ G_{t/n}\big)^n(z)\subseteq B,\quad n\ge 1,
\]
are analytic and bounded uniformly in $t\in\mathbb{S}_{\vartheta_0,\tau}$ and $n\ge 1$ 
(and uniformly on the parameter $0\le\nu\le 1$). 
Since $\mathbb{S}_{\vartheta_0,\tau_2}\subseteq\mathbb{S}_{\vartheta_0,\tau}$ and since the limit \eqref{eq:H-definition}
exists for $t\in\mathbb{S}_{\vartheta_0,\tau_2}$, we conclude from Vitali's theorem (see, e.g., \cite[Proposition 3.2]{BoSi}) 
that for any $z\in O(x_0)\cap Z$ and $0\le\nu\le 1$ the map \eqref{eq:H-curve} extends to an analytic map
\begin{equation}\label{eq:H-curve_extended}
\mathbb{S}_{\vartheta_0,\tau}\to X,\quad t\mapsto H_t^{(\nu)}(z),
\end{equation}
so that 
\begin{equation}\label{eq:H-extended}
\big(S_{t/n}^{(\nu)}\circ G_{t/n}\big)^n(z)\stackrel{X}{\to} H_t^{(\nu)}(z),\quad n\to\infty,
\end{equation}
locally uniformly in $t\in\mathbb{S}_{\vartheta_0,0,\tau}$ and uniformly on $0\le\nu\le 1$.
(Note that, by construction, the range of the map \eqref{eq:H-curve_extended} is contained in 
the closure of the ball $B$ in $X$.)

\begin{Rem}\label{rem:Vitali's_theorem}
One sees from the proof of  \cite[Proposition 3.2]{BoSi} that Vitali's theorem holds for sequences of analytic maps 
that depend on a parameter and converge uniformly in the parameter on an open subset of their domain of definition.
\end{Rem}

\noindent By passing to the limit as $n\to\infty$ in the inequality \eqref{eq:pre-Lipschitz} (with $\alpha_j=\mu_j=1$, 
$1\le j\le \ell=n$), we then conclude that
\begin{equation}\label{eq:Lipschitz_Z}
d\big(H_t^{(\nu)}(z'),H_t^{(\nu)}(z'')\big)\le e^{2\beta|t|}d(z',z'')
\end{equation}
for any $z',z''\in O(x_0)\cap Z$, $t\in\mathbb{S}_{\vartheta_0,0,\tau}$, and $0\le\nu\le 1$.
Now, take $x\in O(x_0)$ and choose a sequence $(z_k)_{k\ge 1}$ such that
$z_k\in O(x_0)\cap Z$ and $z_j\to x$ in $X$ as $j\to\infty$. Then, it follows from \eqref{eq:Lipschitz_Z} that
$d\big(H_t^{(\nu)}(z_j),H_t^{(\nu)}(z_k)\big)\le e^{2\beta|t|}d(z_j,z_k)$
for any $j,k\ge 1$ and for any given $t\in\mathbb{S}_{\vartheta_0,0,\tau}$ and $0\le\nu\le 1$. 
This implies that the sequence $\big(H_t^{(\nu)}(z_j)\big)_{j\ge 1}$ is fundamental in $X$ and hence it converges 
to an element $H_t^{(\nu)}(x)$ in $X$ for any $t\in\mathbb{S}_{\vartheta_0,0,\tau}$ and $0\le\nu\le 1$.
One easily sees that $H_t^{(\nu)}(x)$ is independent of the choice of the sequence $(z_j)_{j\ge 1}$.
By passing to the limit in \eqref{eq:Lipschitz_Z}, we then conclude that
\begin{equation}\label{eq:Lipschitz}
d\big(H_t^{(\nu)}(x),H_t^{(\nu)}(y)\big)\le e^{2\beta|t|}d(x,y)
\end{equation}
for any $x,y\in O(x_0)$, $t\in\mathbb{S}_{\vartheta_0,0,\tau}$, and $0\le\nu\le 1$.
Our next step is to show that for any $x\in O(x_0)$, $0\le\nu\le 1$, and $t\in\mathbb{S}_{\vartheta_0,0,\tau}$, we have that
\begin{equation}\label{eq:H-limit}
\big(S_{t/n}^{(\nu)}\circ G_{t/n}\big)^n(x)\stackrel{X}{\to} H_t^{(\nu)}(x),\quad n\to\infty,
\end{equation}
and that the convergence is uniform in $x$ on compact subsets of $O(x_0)$, 
locally uniform in $t\in\mathbb{S}_{\vartheta_0,0,\tau}$, and uniform on $0\le\nu\le 1$.
To this end, choose $\varepsilon>0$ and let $K\subseteq O(x_0)$ be
a compact set in $X$. For any given $x\in K$ choose $z\in O(x_0)\cap Z$ such that 
$d(x,z)<\rho_\varepsilon$ where $\rho_\varepsilon:=\varepsilon/(8e^{2\beta\tau})$.
In view of \eqref{eq:H-extended}, for any compact set $\mathcal{T}\subseteq\mathbb{S}_{\vartheta_0,\tau}$ we can choose 
$N\ge 1$ such that 
\begin{equation}\label{eq:d<epsilon}
d\big(\big(S_{t/n}^{(\nu)}\circ G_{t/n}\big)^n(z),H_t^{(\nu)}(z)\big)<\varepsilon/2
\end{equation}
for any $n\ge N$, $t\in\mathcal{T}$, and $0\le\nu\le 1$. 
It then follows from \eqref{eq:pre-Lipschitz}, \eqref{eq:Lipschitz}, and \eqref{eq:d<epsilon}, that for any 
$y\in B^d_{\rho_\varepsilon}(x)$ and $n\ge N$ we have that
\begin{align*}
&d\big((S_{t/n}^{(\nu)}\circ G_{t/n}\big)^n(y),H_t^{(\nu)}(y)\big)\le\\
&\le d\big((S_{t/n}^{(\nu)}\circ G_{t/n}\big)^n(y),(S_{t/n}^{(\nu)}\circ G_{t/n}\big)^n(z)\big)
+d\big((S_{t/n}^{(\nu)}\circ G_{t/n}\big)^n(z),H_t^{(\nu)}(z)\big)+
d\big(H_t^{(\nu)}(z),H_t^{(\nu)}(x)\big)\\
&\le e^{2\beta\tau}d\big(y,z\big)+d\big((S_{t/n}^{(\nu)}\circ G_{t/n}\big)^n(z),H_t^{(\nu)}(z)\big)
+e^{2\beta\tau}d\big(z,x\big)<\varepsilon
\end{align*}
for any $t\in\mathcal{T}$ and for any $0\le\nu\le 1$. 
In view of the compactness of the set $K\subseteq O(x_0)$, we then conclude that the convergence in \eqref{eq:H-limit} 
holds uniformly in $x\in K$, $t\in\mathcal{T}$, and $0\le\nu\le 1$, for any choice of the compact subsets
$K\subseteq O(x_0)$ and $\mathcal{T}\subseteq\mathbb{S}_{\vartheta,0,\tau}$. Now, recall that by construction,
for any $n\ge 1$ the map
\[
\mathbb{S}_{\vartheta_0,0,\tau}\times O(x_0)\times[0,1]\to B\subseteq X,\quad
(t,x;\nu)\mapsto(S_{t/n}^{(\nu)}\circ G_{t/n}\big)^n(x),
\]
is continuous, and for any $0\le\nu\le 1$ and $n\ge 1$ the map
\[
\mathbb{S}_{\vartheta_0,\tau}\times O(x_0)\to B\subseteq X,\quad
(t,x)\mapsto(S_{t/n}^{(\nu)}\circ G_{t/n}\big)^n(x),
\]
is analytic. This implies that the flow map
\begin{equation}\label{eq:H-map_in_X}
\mathbb{S}_{\vartheta_0,0,\tau}\times O(x_0)\times[0,1]\to U\subseteq X,\quad
(t,x;\nu)\mapsto H_t^{(\nu)}(x),
\end{equation}
is continuous and that for any $0\le\nu\le 1$ the map
\begin{equation}\label{eq:H-map_in_X_restricted}
\mathbb{S}_{\vartheta_0,\tau}\times O(x_0)\to U\subseteq X,\quad
(t,x;\nu)\mapsto H_t^{(\nu)}(x),
\end{equation}
is analytic. We have the following Lemma.

\begin{Lem}\label{lem:solution}
For any $0\le\nu\le 1$ and $x\in O(x_0)$ the curve
\begin{equation}\label{eq:solution}
[0,\tau)\to X,\quad t\to u^{(\nu)}(t,x):=H_t^{(\nu)}(x),
\end{equation}
is a solution of the equation (cf. \eqref{eq:pde}),
\begin{equation}\label{eq:pde_nu}
\frac{\dd u^{(\nu)}}{\dd t}(t)=\E\big(u^{(\nu)}(t)\big)+\nu\,\Xi\big(u^{(\nu)}(t)\big),\quad u^{(\nu)}\big|_{t=0}=x.
\end{equation}
\end{Lem}

\begin{proof}[Proof of Lemma \ref{lem:solution}]
Take $z_0\in O(x_0)\cap Z$ and let $V_Z(z_0)\subseteq O(x_0)\cap Z$ be the open neighborhood from 
Lemma \ref{lem:commutator_estimate}. We then use Lemma \ref{lem:fundamental} to choose 
open neighborhoods $\widetilde{W}_Z(z_0)\subseteq W_Z(z_0)\subseteq V_Z(z_0)$ of $z_0$ in $Z$
and constants $0<\tau_3\le\tau$ and $K>0$ such that the closure of $W_Z(z_0)$ in $X$ is contained in $O(x_0)\cap Z$
and for any $z\in\widetilde{W}_Z(z_0)$, $t\in\mathbb{S}_{\vartheta_0,0,\tau_3}$, $0\le\nu\le 1$, and $n\ge 1$, we have that
\begin{equation}\label{eq:stability_Z_general'}
\big(S_{t/n}^{(\nu)}\big)^{\alpha_1}\circ G_{t/n}^{\mu_1}\circ\cdots\circ
\big(S_{t/n}^{(\nu)}\big)^{\alpha_\ell}\circ G_{t/n}^{\mu_\ell}(z)\in W_Z(z_0)
\end{equation}
for any multi-indexes $\alpha,\mu\in\Z_{\ge 0}^\ell$ with $1\le\ell\le 2n$, $|\alpha|\le n$, 
$|\mu|\le n$, and
\begin{equation}\label{eq:fundamental'}
d\big(\big(S_{t/n}^{(\nu)}\circ G_{t/n}\big)^n(z),\big(S_{t/m}^{(\nu)}\circ G_{t/m}\big)^m(z)\big)\le 
2 K e^{2\beta|t|}|t|^2/n
\end{equation}
for any $z\in W_Z(z_0)$ and $m\ge n$.
By passing to the limit as $m\to\infty$ in the inequality \eqref{eq:fundamental'}, we obtain that
\begin{equation}\label{eq:fundamental''}
d\big(\big(S_{t/n}^{(\nu)}\circ G_{t/n}\big)^n(z), H^{(\nu)}(z)\big)\le 
2 K e^{2\beta|t|}|t|^2/n
\end{equation}
for any $z\in W_Z(z_0)$, $t\in\mathbb{S}_{\vartheta_0,0,\tau_3}$, $0\le\nu\le 1$, and $n\ge 1$.
We now modify the proof of \cite[Lemma 2.7]{EM}.
For $z\in\widetilde{W}_Z(z_0)$, $t\in\mathbb{S}_{\vartheta_0,0,\tau_3}$, and $0\le\nu\le 1$, we set 
\[
P^{(\nu)}_t(z):=S^{(\nu)}_t\circ G_t(z).
\]
Take $s,t\in[0,\tau_3/2]$. By construction, we then have that for any $z\in\widetilde{W}_Z(z_0)$,
\[
H^{(\nu)}_t\circ H^{(\nu)}_s(z)\quad\text{\rm and}\quad H^{(\nu)}_{s+t}(z)
\]
are well-defined and belong to the ball $B\subseteq X$ where the distance function 
$d : B\times B\to\R_{\ge 0}$ is defined.
Let us first assume that $t$ and $s$ are rationally dependent, i.e., $t=(p/q) s$ where $p,q\in\N$.
Then, for $n\ge 1,$
\[
\frac{s+t}{n}=\frac{s}{n_1}=\frac{t}{n_2}
\]
where $n_1:=\frac{nq}{p+q}$ and $n_2:=\frac{np}{p+q}$. 
In what follows we take $n$ running successively through the infinite set 
\[
A_{p,q}:=\big\{(p+q)l\,\big|\,l\in\N\big\}
\]
of integers. Then,  $n_1=ql$ and $n_2=pl$ (for $l=1,2,...$) are integers and $n_1,n_2\to\infty$ when 
$n\to\infty$, $n\in A_{p,q}$. Since $n_1+n_2=n$, we have
\begin{equation}\label{eq:eqlt1}
(P^{(\nu)}_{\frac{s+t}{n}})^n(z)=(P^{(\nu)}_{\frac{s+t}{n}})^{n_1}\circ(P^{(\nu)}_{\frac{s+t}{n}})^{n_2}(z)=
(P^{(\nu)}_{s/n_1})^{n_1}\circ(P^{(\nu)}_{t/n_2})^{n_2}(z)\,.
\end{equation}
For any $z\in\widetilde{W}_Z(z_0)$, $t\in[0,\tau_3/2]$, and for any $n\in A_{p,q}$, $n_1$, 
and $n_2$ as above, we have
\begin{align*}
&d\big((P^{(\nu)}_{s/n_1})^{n_1}\circ(P^{(\nu)}_{t/n_2})^{n_2}(z),H^{(\nu)}_s\circ H^{(\nu)}_t(z)\big)\le
d\big((P^{(\nu)}_{s/n_1})^{n_1}\circ(P^{(\nu)}_{t/n_2})^{n_2}(z),H^{(\nu)}_s\circ(P^{(\nu)}_{t/n_2})^{n_2}(z)\big)\\
&+d\big(H^{(\nu)}_s\circ(P^{(\nu)}_{t/n_2})^{n_2}(z),H^{(\nu)}_s\circ H^{(\nu)}_t(z)\big)\le\frac{C_1}{n_1}+
e^{2\beta\tau_1}\,d\big((P^{(\nu)}_{t/n_2})^{n_2}(z),H^{(\nu)}_t(z)\big)=C_2\Big(\frac{1}{n_1}+\frac{1}{n_2}\Big)
\end{align*}
where we used  \eqref{eq:Lipschitz}, \eqref{eq:stability_Z_general'}, and \eqref{eq:fundamental''},
and the constants are independent of the choice of $z\in\widetilde{W}_Z(z_0)$ and $t\in[0,\tau_3/2]$.
Hence, $(P^{(\nu)}_{s/n_1})^{n_1}\circ(P^{(\nu)}_{t/n_2})^{n_2}(z)$ converges to 
$H^{(\nu)}_s\circ H^{(\nu)}_t(z)$ in $X$ as $n\to\infty$ with $n\in A_{p,q}$. 
By passing to the limit as $n\to\infty$ with $n\in A_{p,q}$ in \eqref{eq:eqlt1}, we then conclude 
from \eqref{eq:H-limit} (with $t$ replaced by $s+t$) that $H^{(\nu)}_{s+t}(z)=H^{(\nu)}_s\circ H^{(\nu)}_t(z)$ in 
the case when $s$ and $t$ are rationally dependent.
Since the rational numbers are dense in $\R$, we obtain from the continuity
of the map \eqref{eq:H-map_in_X} that
\begin{equation}\label{eq:H-the_semigroup_property}
H^{(\nu)}_{s+t}(z)=H^{(\nu)}_s\circ H^{(\nu)}_t(z)
\end{equation}
for any $z\in\widetilde{W}_Z(z)$, $s,t\in[0,\tau_3/2)$, and $0\le\nu\le 1$.
It follows from \eqref{eq:fundamental''} that for any $z\in\widetilde{W}_Z(z_0)$ we have that
\[
d\big(S_t^{(\nu)}\circ G_t(z), H^{(\nu)}(z)\big)=O(|t|^2)
\]
uniformly in $t\in[0,\tau_3)$ and $0\le\nu\le 1$. This implies that
\[
\frac{\dd}{\dd t}\Big|_{t=0}\big[H^{(\nu)}_t(z)\big]=\frac{\dd}{\dd t}\Big|_{t=0}\big[S_t^{(\nu)}\circ G_t(z)\big]
=\E(z)+\nu\,\Xi(z)
\]
for any $z\in\widetilde{W}_Z(z_0)$ (cf. Corollary \ref{coro:C^k_along_curves}). 
This together with \eqref{eq:H-the_semigroup_property} implies that
\begin{equation}\label{eq:pde_local}
\frac{\dd}{\dd t}\big[H^{(\nu)}_t(z)\big]=\E\big(H^{(\nu)}_t(z)\big)+\nu\,\Xi\big(H^{(\nu)}_t(z)\big)
\end{equation}
for any $z\in\widetilde{W}_Z(z_0)$, $t\in(0,\tau_3/2)$, and $0\le\nu\le 1$.
In view of the analyticity of the maps \eqref{eq:H-map_in_X_restricted}, $\Xi : U\to\widetilde{X}$, and $\E : U\to X$, 
we then conclude from \eqref{eq:pde_local}, the continuity of the map \eqref{eq:H-map_in_X}, 
and the fact that $Z$ is dense in $X$, that \eqref{eq:pde_local} holds for any $x\in O(x_0)$, 
$t\in[0,\tau)$, and $0\le\nu\le 1$.
This completes the proof of the Lemma.
\end{proof}

By the analyticity in $t\in\mathbb{S}_{\vartheta_0,\tau}$ of the both sides of \eqref{eq:pde_nu}, 
the equation holds for any $x\in O(x_0)$ and $t\in\mathbb{S}_{\vartheta_0,\tau}$. 
Since $u^{(\nu)}\in C_b\big(\mathcal{S}_{\vartheta_0,0,\tau},U\big)$, we then conclude that
\[
u^{(\nu)}(t)=x+\int_0^t\Big(\E\big(u^{(\nu)}(s)\big)+\nu\,\Xi\big(u^{(\nu)}(s)\Big)\,ds
\]
where the integral exists as a Riemann integral in $\widetilde{X}$ and the integration is 
taken over the straight segment $[0,t]\subseteq\mathcal{S}_{\vartheta_0,0,\tau}$. In particular, this implies that
\[
\frac{\dd u^{(\nu)}}{\dd t}(0)=\lim_{t\to 0, t\in\mathcal{S}_{\vartheta_0,\tau}}
\frac{u^{(\nu)}(t)-x}{t}
\]
where the limit is taken in $\widetilde{X}$.
This completes the proof of Proposition \ref{prop:product_formula}.
\end{proof}

\begin{proof}[Proof of Theorem \ref{th:product_formula}]
The Theorem follows from Proposition \ref{prop:product_formula} and Lemma \ref{lem:equivalent_distance}. 
In fact, the estimates \eqref{eq:E-estimates} in Lemma \ref{lem:equivalent_distance} hold locally in 
$U$ (resp. $U_Z$) since the maps in \eqref{eq:G} are assumed $C^2$ on conic sets. 
The global in time estimates \eqref{eq:dS_t-exponential_growth1} and \eqref{eq:dS_t-exponential_growth2} 
then allow us to apply Lemma \ref{lem:equivalent_distance} twice (with $N\equiv N_X(x_0)$, and 
$N\equiv N_Z(z_0)$, respectively) and conclude that the quasi-contraction condition $(b)$ in 
Proposition \ref{prop:product_formula} holds.
The uniqueness of the solution \eqref{eq:solution} follows easily from the arguments in 
the proof of \cite[Lemma 6]{EM}. Theorem \ref{th:product_formula} then follows from 
Proposition \ref{prop:product_formula}.
\end{proof}

\appendix

\section{The class of $C^k$ maps on conic sets}\label{sec:C^k_along_curves}
In this Appendix we define and study the main properties of the maps of class {\em $C^k$ on conic sets} 
in the context of complex Banach spaces. 
Take $0<T\le\infty$ and consider the conic set in $\C$ (cf. \eqref{eq:other_cones}),
\begin{equation}\label{eq:the_T-cone}
\mathbb{S}_{\vartheta,0,T}=\big\{z\in\mathbb{S}_{\vartheta,0}\,\big|\,|z|<T\big\}.
\end{equation}
Let $F : \mathbb{S}_{\vartheta,0,T}\times U\to Y$, where $X$ and $Y$ are complex Banach spaces and $U$ is a non-empty
open set in $X$, be a continuous map. For a given $x\in U$, the partial derivative of $\big(D_1F\big)(t,x)$ in the direction of 
$t\in\mathbb{S}_{\vartheta}$ at $t=0$ is defined as
\begin{equation}\label{eq:derivative_at_t=0}
\big(D_1F\big)(0,x):=
\lim_{s\to 0, s\in\mathbb{S}_{\vartheta,T}}\frac{F(s,x)-F(0,x)}{s}
\end{equation}
provided that the limit exists in $Y$ for $s\to 0$ with $s\in\mathbb{S}_{\vartheta,0,T}$.
The derivative $\big(D_1F\big)(t,x)$ for $t\in\mathbb{S}_{\vartheta,T}$ is defined in a standard way.
The directional partial derivative of $F$ at $(t,x)\in\mathbb{S}_{\vartheta,0,T}\times U$ in the direction of 
$\xi\in X$ is defined as
\[
\big(D_2f\big)(t,x)\big(\xi\big):=
\lim_{s\to 0, s\in\C}\frac{F(t,x+s\xi)-F(t,x)}{s}
\]
provided that the limit exists in $Y$. Assume that the limit above exists for all 
$(t,x)\in\mathbb{S}_{\vartheta,0,T}\times U$ and $\xi\in X$.
Then the second order directional partial derivative of $F$ at $(t,x)\in\mathbb{S}_{\vartheta,0,T}\times U$ in 
the direction of $\xi_1,\xi_2\in X$ is defined as
\[
\big(D_2^2f\big)(t,x)\big(\xi_1,\xi_2\big):=\lim_{s\to 0, s\in\C}
\frac{\big(D_2F\big)(t,x+s\xi_2)(\xi_1)-\big(D_2F\big)(t,x)(\xi_1)}{s}
\]
provided that the limit exists in $Y$. Similarly, the second order directional partial derivative 
$(D_1D_2F)(0,x)(\xi)$ is defined as in \eqref{eq:derivative_at_t=0},
\[
\big(D_1D_2F\big)(0,x)\big(\xi\big):=\lim_{s\to 0, s\in\mathbb{S}_{\vartheta,T}}
\frac{\big(D_2F\big)(s,x)\big(\xi\big)-\big(D_2F\big)(0,x)\big(\xi\big)}{s}\,.
\]
For $t\in\mathbb{S}_{\vartheta,T}$ the derivative $\big(D_1D_2F\big)(t,x)\big(\xi\big)$ is defined in a standard way.
The higher order directional partial derivatives of $F$ at $(t,x)\in\mathbb{S}_{\vartheta,0,T}\times U$ are 
then defined recursively.

\begin{Def}\label{def:C^k_along_curves_complex}
A continuous map $F : \mathbb{S}_{\vartheta,0,T}\times U\to Y$ is {\em $C^k$ on conic sets}
for some $k\ge 0$ if for any multi-indexes $\mu,\alpha\in\Z_{\ge 0}^\ell$, $|\mu|+|\alpha|\le k$, $\ell\in\Z_{\ge 0}$, 
and directions $\xi_1,...,\xi_{|\alpha|}$, there exists a directional partial derivative
\[
\big(D^{\mu,\alpha}F\big)(t,x)\big(\xi_1,...,\xi_{|\alpha|}\big):=
\big(D_1^{\mu_1}D_2^{\alpha_1}\ldots D_1^{\mu_\ell}D_2^{\alpha_\ell}F\big)(t,x)\big(\xi_1,...,\xi_{|\alpha|}\big)
\]
and the map
\begin{equation}\label{eq:D_1D_2-continuity}
\mathbb{S}_{\vartheta,0,T}\times U\times X^{|\alpha|}\to Y,\quad
(t,x,\xi)\mapsto\big(D^{\mu,\alpha}F\big)(t,x)\big(\xi_1,...,\xi_{|\alpha|}\big),
\end{equation}
is continuous.
\end{Def}

In the same way as above we define the derivative $\frac{\dd\gamma}{\dd t}(t)$ of a map
$\gamma : \mathcal{S}_{\vartheta,0,T}\to Y$.

\begin{Def}\label{def:C^k_complex}
A continuous map $\gamma : \mathcal{S}_{\vartheta,0,T}\to Y$ is $C^k$ for some $k\ge 0$ if 
for any $0\le j\le k$ and $t\in\mathbb{S}_{\vartheta,0,T}$ the derivative $\frac{\dd^j\gamma}{\dd t^j}(t)$ exists 
and $\frac{\dd^j\gamma}{\dd t^j} : \mathbb{S}_{\vartheta,0,T}\to Y$ is continuous.
\end{Def}

\noindent The maps of class $C^k$ on conic sets will be also called {\em $C^k$ maps on conic sets}.
By the definition above, a map is $C^0$ on conic sets if and only if it is continuous. 
We will see below that $C^k$ maps on conics sets with $k\ge 2$ have commuting directional partial derivatives, i.e.,
\begin{equation}\label{eq:commuting_derivatives}
\big(D_1^{\mu_1}D_2^{\alpha_1}\ldots D_1^{\mu_\ell}D_2^{\alpha_\ell}F\big)(t,x)\big(\xi_1,...,\xi_{|\alpha|}\big)=
\big(D_1^{|\mu|}D_2^{|\alpha|}F\big)(t,x)\big(\xi_1,...,\xi_{|\alpha|}\big)
\end{equation}
where the derivative on the right side is symmetric in $\xi_1,...,\xi_{|\alpha|}$.

\begin{Rem}\label{rem:C^1=>analytic}
It follows from Definition \ref{def:C^k_along_curves_complex} that $C^k$ maps on conic sets with $k\ge 1$
are holomorphic on affine lines and locally uniformly bounded on $\mathbb{S}_{\vartheta,T}\times U$. 
Hence, any such map $F : \mathbb{S}_{\vartheta,0,T}\times U\to X$ is (complex) analytic in the interior
$\mathbb{S}_{\vartheta,T}\times U$ of the domain (see, e.g., \cite[Theorem 1, Appendic A]{PT}). 
In particular, for any given $t\in\mathbb{S}_{\vartheta,T}$ and $x\in U$ the directional partial derivative 
$\big(D_2F\big)(t,x)\big(\xi\big)$ in the direction of $\xi\in X$ coincides with the Fr\'echet derivative 
$\big(\dd_x F(t,\cdot)\big)(\xi)$.
\end{Rem}

\noindent Definition \ref{def:C^k_along_curves_complex} means simply that the directional partial derivatives 
of $F$ of order $\le k$ (with any order of appearance of $D_1$ and $D_2$) exist 
and are weakly continuous on $\mathbb{S}_{\vartheta,0,T}\times U$. The following Lemma gives an alternative 
characterization of the maps of class $C^k$ on conic sets with $k\ge 1$.

\begin{Lem}\label{lem:C^k_along_curves}
A continuous map $F : \mathbb{S}_{\vartheta,0,T}\times U\to Y$ is {\em $C^k$ on conic sets} for some $k\ge 1$ 
if and only if its restriction $F_0:=F\big|_{\mathbb{S}_{\vartheta,T}\times U}$ satisfies the following properties:
\begin{itemize}
\item[$(i)$] The map $F_0 : \mathbb{S}_{\vartheta,T}\times U\to Y$ is analytic.

\item[$(ii)$] For any multi-index $\alpha\in\Z^2_{\ge 0}$, $0\le|\alpha|\le k$, 
and for any $\xi_1,...,\xi_{\alpha_2}\in X$ the derivative
\[
\mathbb{S}_{\vartheta,T}\times U\times X^{\alpha_2}\to Y,\quad
(t,x,\xi)\mapsto\big(D_1^{\alpha_1}D_2^{\alpha_2}F_0\big)(t,x)\big(\xi_1,...,\xi_{\alpha_2}\big)
\]
extends to a continuous map on $\mathbb{S}_{\vartheta,0,T}\times U$.
\end{itemize}
If a map $F_0 : \mathbb{S}_{\vartheta,T}\times U\to Y$ satisfies conditions $(i)$ and $(ii)$ then
its continuous extension $F : \mathbb{S}_{\vartheta,0,T}\times U\to Y$ is $C^k$ on conic sets.
\end{Lem}

\begin{proof}[Proof of Lemma \ref{lem:C^k_along_curves}]
The direct statement follows from Remark \ref{rem:C^1=>analytic}. Let us prove the converse one. 
Assume that for some $k\ge 1$ the map $F_0 : \mathbb{S}_{\vartheta,T}\times U\to Y$ satisfies 
conditions $(i)$ and $(ii)$ of the Lemma. Denote by $F : \mathbb{S}_{\vartheta,0,T}\times U\to Y$ its continuous 
extension on $\mathbb{S}_{\vartheta,0,T}\times U$ and take $0\le\alpha_2\le k-1$, $\alpha_2\in\Z_{\ge 0}$. 
Then, by $(i)$, for any given $x\in U$ and $\xi_1,...,\xi_{\alpha_2}\in X$ the map
\begin{equation}\label{eq:z-map}
\mathbb{S}_{\vartheta,T}\to Y,\quad t\mapsto z_{\alpha_2}(t):=\big(D_2^{\alpha_2} F_0\big)(t,x)\big(\xi_1,...,\xi_{\alpha_2}\big),
\end{equation}
is holomorphic. Hence, for any given $s,t\in\mathbb{S}_{\vartheta,T}$ we have that
\begin{equation}\label{eq:t-derivative1}
z_{\alpha_2}(t)=z_{\alpha_2}(s)+\int_s^t\frac{\dd z_{\alpha_2}}{\dd t}(\lambda)\,d\lambda
\end{equation}
where $\frac{\dd z_{\alpha_2}}{\dd t}(\lambda)=\big(D_1D_2^{\alpha_2} F_0\big)(\lambda,x)\big(\xi_1,...,\xi_{\alpha_2}\big)$
and the integration is taken over the straight segment $[s,t]$ in $\C$. It then follows from $(ii)$ that we can pass to the limit
in \eqref{eq:t-derivative1} as $s\to 0$, $s\in\mathbb{S}_{\vartheta,T}$, to conclude that
$z_{\alpha_2}(t)=z_{\alpha_2}(0)+\int_0^t\frac{\dd z_{\alpha_2}}{\dd t}(\lambda)\,d\lambda$ where
\begin{equation}\label{eq:z-limit1}
z_{\alpha_2}(0):=\lim_{s\to 0,s\in\mathbb{S}_{\vartheta,T}}z_{\alpha_2}(s)\,.
\end{equation}
This implies that
\begin{equation}\label{eq:t-derivative2}
\frac{z_{\alpha_2}(t)-z_{\alpha_2}(0)}{t}-\frac{\dd z_{\alpha_2}}{\dd t}(0)=\frac{1}{t}
\int_0^t\Big(\frac{\dd z_{\alpha_2}}{\dd t}(\lambda)-\frac{\dd z_{\alpha_2}}{\dd t}(0)\Big)\,d\lambda
\end{equation}
where
\begin{equation}\label{eq:z-limit2}
\frac{\dd z_{\alpha_2}}{\dd t}(0):=\lim_{s\to 0,s\in\mathbb{S}_{\vartheta,T}}\frac{\dd z_{\alpha_2}}{\dd t}(s)\,.
\end{equation}
It then follows from \eqref{eq:z-limit2} that we can pass to the limit in \eqref{eq:t-derivative2} as 
$t\to 0$, $t\in\mathbb{S}_{\vartheta,T}$, to conclude that
\begin{equation}\label{eq:z-limit2=t-derivative}
\frac{\dd z_{\alpha_2}}{\dd t}(0)=\lim_{t\to 0,t\in\mathbb{S}_{\vartheta,T}}\frac{z_{\alpha_2}(t)-z_{\alpha_2}(0)}{t}\,.
\end{equation}
We will prove below that for any given integer $\alpha_2\ge 0$,
$x\in U$, and $\xi_1,...,\xi_{\alpha_2}\in X$, the directional partial derivative 
$\big(D_2^{\alpha_2} F\big)(0,x)\big(\xi_1,...,\xi_{\alpha_2}\big)$ exists and that
\begin{equation}\label{eq:z_limit1=derivative}
\big(D_2^{\alpha_2} F\big)(0,x)\big(\xi_1,...,\xi_{\alpha_2}\big)=
\lim_{t\to 0,t\in\mathbb{S}_{\vartheta,T}}\big(D_2^{\alpha_2} F_0\big)(t,x)\big(\xi_1,...,\xi_{\alpha_2}\big)\,.
\end{equation}
It then follows from \eqref{eq:z-map} and \eqref{eq:z-limit2=t-derivative} that
$\big(D_1D_2^{\alpha_2} F\big)(0,x)\big(\xi_1,...,\xi_{\alpha_2}\big)$ exists for $0\le\alpha_2\le k-1$ and that
\begin{equation}\label{eq:z_limit2=derivative}
\big(D_1D_2^{\alpha_2} F\big)(0,x)\big(\xi_1,...,\xi_{\alpha_2}\big)=
\lim_{t\to 0,t\in\mathbb{S}_{\vartheta,T}}\big(D_1D_2^{\alpha_2} F_0\big)(t,x)\big(\xi_1,...,\xi_{\alpha_2}\big)\,.
\end{equation}
Let us now prove that \eqref{eq:z_limit1=derivative} holds. 
Since $F : \mathbb{S}_{\vartheta,0,T}\times U\to Y$ is continuous, we conclude that
\[
F_0(t,\cdot)\to F(0,\cdot)\quad\text{\rm as}\quad t\to 0,\,\,t\in\mathbb{S}_{\vartheta,0,T},
\]
where the convergence in $Y$ is uniform in $x\in K$ for any compact subset $K\subseteq U$.
Since for any given $t\in\mathbb{S}_{\vartheta,T}$ the function $F_0(t,\cdot) : U\to Y$ is analytic, 
we then obtain from Weierstrass theorem that the limit function $F(0,\cdot) : U\to Y$ is analytic on affine lines,
and hence it is analytic. This also shows that for any $x\in U$ and $\xi\in X$ we have that
\[
\frac{\dd^jF}{\dd s^j}(0,x+s\xi)\Big|_{s=0} =
\lim_{t\to 0,t\in\mathbb{S}_{\vartheta,T}}\Big(\frac{\dd^j F_0}{\dd s^j}(t,x+s\xi)\Big|_{s=0}\Big)
\]
for any $j\ge 0$. This implies that \eqref{eq:z_limit1=derivative} holds with $\xi_1=...=\xi_{\alpha_2}=\xi$
for any $\alpha_2\in\Z_{\ge 0}$ and $x\in U$. Since the functions $F(0,\cdot) : U\to Y$ and
$F_0(t,\cdot) : U\to Y$, $t\in\mathbb{S}_{\vartheta,T}$, are analytic, the multilinear forms in
\eqref{eq:z_limit1=derivative} are symmetric. This proves that \eqref{eq:z_limit1=derivative} holds.
In view of the continuity in $(ii)$, the limits in \eqref{eq:z_limit1=derivative} and \eqref{eq:z_limit2=derivative}
exist not only as $t\to 0$, $t\in\mathbb{S}_{\vartheta,T}$, for fixed $x\in U$ but also as
$(t,x)\to(0,x_0)$ in $\mathbb{S}_{\vartheta}\times U$, $x_0\in U$.
By iterating these arguments we prove that $F$ is $C^k$ on conic sets.
\end{proof}

Let $Z$ be a complex Banach space and $V\subseteq Y$ be an open set in $Y$. 

\begin{Coro}\label{coro:C^k_along_curves}
Assume that the maps $F : \mathbb{S}_{\vartheta,0,T}\times U\to V$ and $G : \mathbb{S}_{\vartheta,0,T}\times V\to Z$
are $C^k$ on conic sets. Then the composed map $\mathbb{S}_{\vartheta,0,T}\times U\to Z$, 
$(t,x)\mapsto G\big(t,F(t,x)\big)$, is $C^k$ on conic sets and its derivatives of order $\le k$ satisfy the product
rule on $\mathbb{S}_{\vartheta,0,T}\times U$.
\end{Coro}

\begin{proof}[Proof of Corollary \ref{coro:C^k_along_curves}]
For $(t,x)\in\mathbb{S}_{\vartheta,0,T}\times U$ denote $(G\circ F)(t,y):=G\big(t,F(t,x)\big)$.
Then, $G\circ F\big|_{\mathbb{S}_{\vartheta,T}\times U}=G_0\circ F_0$.
Since the composition of analytic maps is an analytic map, we then conclude that 
the restricted map $G\circ F\big|_{\mathbb{S}_{\vartheta,T}\times U}$ is analytic, and hence it satisfies condition $(i)$. 
It follows from the product rule and the fact that $F_0$ and $G_0$ satisfy condition $(ii)$ that 
the restriction $G\circ F\big|_{\mathbb{S}_{\vartheta,T}\times U}$ also satisfies condition $(ii)$. 
The Corollary then follows from Lemma \ref{lem:C^k_along_curves}.
\end{proof}

The following Corollary shows that the directional partial derivatives $D_1$ and
$D_2$ of a $C^k$ on conic sets with $k\ge 2$ commute not only in the interior of 
the domain but also at the points on the boundary (cf. \eqref{eq:commuting_derivatives}).

\begin{Coro}\label{coro:commuting_derivatives}
If $F : \mathbb{S}_{\vartheta,0,T}\times U\to Y$ is $C^k$ on conic sets for some $k\ge 2$ then
the directional partial derivatives $D_1$ and $D_2$ of $F$ commute.
\end{Coro}

\noindent The Corollary follows from Lemma \ref{lem:C^k_along_curves}$(i)$ and the continuity of 
the map \eqref{eq:D_1D_2-continuity}.
We will also need the following Proposition.

\begin{Prop}\label{prop:analytic_derivatives}
If $F : \mathbb{S}_{\vartheta,0,T}\times U\to Y$ is $C^k$ on conic sets for some $k\ge 1$ then
the derivatives $\big(D_1^j F\big)(0,\cdot) : U\to Y$ are analytic for $0\le j\le k$.
\end{Prop}

\begin{proof}[Proof of Proposition \ref{prop:analytic_derivatives}]
Recall from Lemma \ref{lem:C^k_along_curves} that the restricted map 
$F_0=F\big|_{\mathbb{S}_{\vartheta,T}\times U}$ is analytic.
The analyticity of the $t$ derivatives at $t=0$ then follows since, from the continuity of the map \eqref{eq:D_1D_2-continuity},
for any $0\le j\le k$ we have that $\big(D_1^j F\big)(t,\cdot)\to\big(D_1^j F\big)(0,\cdot)$ as $t\to 0$, 
$t\in\mathbb{S}_{\vartheta,T}$, where the converges in $Y$ is uniformly on compact subsets of $U$, 
and since for any given $t\in\mathbb{S}_{\vartheta,T}$ the derivatives 
$\big(D_1^j F\big)(t,\cdot) : U\to Y$ are analytic by the analyticity of the restriction $F_0$.
\end{proof}

\medskip

Definition \ref{def:C^k_along_curves_complex} is a holomorphic version of \cite[Definition 1]{McOwenTopalov6} in 
the case of real Banach spaces. In that case, Lemma \ref{lem:C^k_along_curves} continues to hold with the condition $(i)$ 
replaced by the assumption that $F_0 : (0,T)\times U\to Y$ is a $C^k$ map.
As in Corollary \ref{coro:C^k_along_curves}, one sees that the composition of $C^k$ maps on conic sets is 
again a $C^k$ map on conic sets (cf., e.g., \cite[Theorem 3.6.4]{Hamilton}).


\begin{thebibliography}{99}   

\bibitem{Arnold} V. Arnold, {\em Sur la geometri{\'e} differentielle des groupes de Lie de dimension infinie et 
ses applications {\`a} l'hydrodynamique des fluids parfaits}, Ann. Inst. Fourier, $\bf 16$(1966), no. 1, 319-361

\bibitem{BoSi} J. Bochnak, J. Siciak, {\em Analytic functions in topological vector spaces},
Studia Mathematica, $\bf 39$(1971),  77-112

\bibitem{BrezisPazy} H. Brezis, A. Pazy, {Convergence and approximation of nonlinear 
operators of semigroups in Banach Spaces},  J. Funct. Anal., $\bf 9$(1972), 63-74

\bibitem{EM} D. Ebin, J. Marsden, {\em Groups of diffeomorphisms and the motion of an incompressible fluid}, 
Ann. Math., $\bf 92$(1970), 102-163

\bibitem{Gunther}  N. G\"unther, {\em \"Uber ein Hauptproblem der Hydrodynamik},
Math. Z., $\bf 24$(1926), no. 1, 448–499



\bibitem{Hamilton} R. Hamilton, {\em The Inverse Function Theorem of Nash and Moser},
Bulletin of the AMS, $\bf 7$(1982), no. 1,  65-222

\bibitem{HaydenMassey} T. Hayden, F. Massey, {\em Nonlinear holomorphic semigroups},
Pacific J. Math., $\bf 57$(1975), no. 2, 423-439

\bibitem{IKT} H. Inci, T. Kappeler and P. Topalov, {\it On the regularity of the composition of 
diffeomorphisms}, \emph{Mem. Amer. Math. Soc.}, \textbf{226} (2013), no.\ 1062


\bibitem{Kato0} T. Kato, {\em Nonlinear semigroups and evolution equations},
J. Math. Soc. Japan, $\bf 19$(1967), no. 4, 508-519

\bibitem{Kato1} T. Kato, {\em Nonstationary flows of viscous and ideal fluids in $\R^3$},
J. Func. Anal., $\bf 9$ (1972), 296-305
 
\bibitem{Kato2} T. Kato, {\em Remarks on zero viscosity limit for nonstationary Navier-Stokes flows with boundary},
Seminar on Nonlinear Partial Differential Equations, Springer, 1984, 85-98

\bibitem{Komura} Y. Komura, {\em Nonlinear semi-groups in Hilbert space}, J. Math. Soc. Japan,
$\bf 19$(1967), no. 4, 493-507

\bibitem{LL} G. Ladas, V. Lakshmikantham, Differential Equations in Abstract Spaces,
Mathematics in Science and Engineering, $\bf 85$, Academic Press, 1972

\bibitem{Lang} S. Lang, {Fundamentals in Differential Geometry}, Graduate Texts in Mathematics, $\bf 191$,
Springer, 1999

\bibitem{McGrath} F. McGrath, {\em Nonstationary plane flow of viscous and ideal fluid},
Arch. Rational Mech. Anal., $\bf 27$(1967), 329–348
 
\bibitem{MB} A. Majda, A. Bertozzi, Vorticity and Incompressible Flow,
Cambridge Texts in Applied Mathematics, Cambridge University Press, 2002

\bibitem{Marsden} J. Marsden, {\em On the product formulas for nonlinear semigroups},
J. Funct. Anal., $\bf 13$(1974), 51-72

\bibitem{McOwen} R. McOwen, {\em The behavior of the Laplacian on weighted Sobolev spaces},
Comm. Pure Appl. Math., $\bf 32$(1979), no. 6, 783-795

\bibitem{McOwenTopalov2} R. McOwen, P. Topalov, 
{\em Groups of asymptotic diffeomorphisms}, Discrete Contin. Dyn. Syst., $\bf 36$ (2016), no. 11, 6331-6377

\bibitem{McOwenTopalov3} R. McOwen, P. Topalov, 
{\em Spatial asymptotic expansions in the incompressible Euler equation}, 
GAFA, $\bf 27$ (2017), 637-675

\bibitem{McOwenTopalov4} R. McOwen, P. Topalov,
{\em Perfect fluid flows on $\R^d$ with growth/decay conditions at infinity},
Mathematische Annalen, $\bf 383$(2022),  1451-1488

\bibitem{McOwenTopalov5} R. McOwen, P. Topalov,
{\em Spatial asymptotic expansions in the Navier-Stokes equation},
IMRN 2024, no. 4, 3391-3441

\bibitem{McOwenTopalov6} R. McOwen, P. Topalov,
{\em The viscosity limit of fluid flows with growth/decay conditions at infinity},
arXiv:2412.05715

\bibitem{MY} G. Misiolek, T. Yoneda, {\em Continuity of the solution map of the Euler equation 
in H{\"o}lder spaces and weak norm inflation in Besov spaces}, Trans. Amer. Math. Soc.,
$\bf 370$(2018), no. 7, 4709-4730

\bibitem{Nelson} E. Nelson, {\em Feynman integrals and the Schr\"odinger equation},
J. Math. Phys., $\bf 5$(1964), no. 3,  332-343

\bibitem{Oharu1} S. Oharu, {\em On the generation of semigroups of nonlinear contractions},
J. Math. Soc. Japan, $\bf 22$(1970), no. 4, 526-550

\bibitem{Pazy} A. Pazy, Semigroups of Linear Operators and Applications to Partial Differential 
Equations, Applied Mathematical Sciences, $\bf 44$, Springer-Verlag, 1983

\bibitem{PT} J. P\"oschel, E.Trubowitz, Inverse Spectral Theory, Academic Press, 1987


\bibitem{Swann} H. Swann, {\em The convergence with vanishing viscosity of nonstationary Navier-Stokes flow
to ideal flow in $\R_3$}, Transactions of the AMS, $\bf 157$(1971), 373-397 

\bibitem{Takhtajan} L. Takhtajan, Quantum Mechanics for Mathematicians, Graduate Studies in Mathematics,
$\bf 95$, AMS, 2008

\bibitem{Top1} P. Topalov, {\em Spatial decay/asymptotics in the Navier-Stokes equation}, 
 Russ. J. Math. Phys., $\bf 32$(2025), no. 1, 196-209
 
 \bibitem{Trotter} H. Trotter, {\em On the product of semi-groups of operators},
 Proc. Amer. Math. Soc., $\bf 10$(1959), 545–551


\bibitem{Wol} W. Wolibner, {\em Un theor\`eme sur l'existence du mouvement plan d'un fluide parfait, homog\`ene, 
incompressible, pendant un temps infiniment long}, Math. Z., $\bf 37$(1933), no. 1, 698-726


\end{thebibliography}
\end{document}